\documentclass{gtmon_a}
\pdfoutput=1

\usepackage{epsfig}
\usepackage{float}
\usepackage{amsmath}
\usepackage{latexsym}
\usepackage{amssymb}
\usepackage{pinlabel}
\usepackage[all]{xy}
\input{diagxy}

\proceedingstitle{The interaction of finite-type and Gromov--Witten
invariants (BIRS 2003)}
\conferencestart{15 November 2003}
\conferenceend{20 November 2003}
\conferencename{The interaction of finite-type and Gromov--Witten
invariants}
\conferencelocation{Banff International Research Station, Banff, Alberta,
Canada}

\editor{David Auckly}
\givenname{David}
\surname{Auckly}

\editor{Jim Bryan}
\givenname{Jim}
\surname{Bryan}

\title{Hodge-type integrals on moduli spaces of admissible covers}

\author{Renzo Cavalieri}
\givenname{Renzo}
\surname{Cavalieri}
\address{Department of Mathematics\\
University of Utah\\\newline
155 South 1400 East\\
Salt Lake City UT 84112\\
USA}
\address{Department of Mathematics\\
University of Michigan\\\newline
2074 East Hall\\
530 Church Street\\
Ann Arbor MI 48109-1043\\
USA}
\email{crenzo@umich.edu}
\email{renzo@math.utah.edu}
\urladdr{}

\volumenumber{8}
\issuenumber{}
\publicationyear{2006}
\papernumber{8}
\startpage{167}
\endpage{194}

\doi{}
\MR{}
\Zbl{}

\arxivreference{math.AG/0411500}


\received{9 March 2005}
\revised{}
\accepted{10 March 2005}
\published{21 September 2007}
\publishedonline{21 September 2007}
\proposed{}
\seconded{}
\corresponding{}
\editor{}
\version{}

\newcommand{\proj}{\mathbb{P}^1}

\newcommand{\Admunz}{\overline{\Adm}({0\stackrel{3}{\longrightarrow}0,((3),t_1,t_2)})}
\newcommand{\Admun}{\overline{\Adm}({g\stackrel{3}{\longrightarrow}0,((3),t_1,\ldots,t_{2g+2})})}
\newcommand{\Admunt}{\overline{\Adm}({g\stackrel{3}{\longrightarrow}0,(t_1,\ldots,t_{2g+4})})}
\newcommand{\Admuntz}{\overline{\Adm}({0\stackrel{3}{\longrightarrow}0,(t_1,\ldots,t_{4})})}
\newcommand{\Admungone}{\overline{\Adm}({g_1\stackrel{3}{\longrightarrow}0,((3),t_1,\ldots,t_{2g_1+2})})}
\newcommand{\Admungtwo}{\overline{\Adm}({g_2\stackrel{3}{\longrightarrow}0,((3),t_1,\ldots,t_{2g_2+2})})}
\newcommand{\Admuntwogtwo}{\overline{\Adm}({g_2\stackrel{2}{\longrightarrow}0,(t_1,\ldots,t_{2g_2+2})})}
\newcommand{\Admuntwog}{\overline{\Adm}({g\stackrel{2}{\longrightarrow}0,(t_1,\ldots,t_{2g+2})})}
\newcommand{\Admuntgone}{\overline{\Adm}({g_1\stackrel{3}{\longrightarrow}0,(t_1,\ldots,t_{2g_2+4})})}

\newcommand{\E}{R^1\pi_\ast f^\ast(\mathcal{O}_{\proj}\oplus\mathcal{O}_{\proj}(-1))}
\newcommand{\F}{R^1\pi_\ast f^\ast(\mathcal{O}_{\proj}(-1)\oplus\mathcal{O}_{\proj}(-1))}
\newcommand{\Admtwog}{\overline{\Adm}
  ({g\stackrel{2}{\longrightarrow}\mathbb{P}^1_{\mathbb{C}},
  (t_1,t_2,\ldots,t_{2g+2})})}

\newcommand{\Admtwoun}{\overline{\Adm}
  ({g\stackrel{2}{\longrightarrow}0,(t_1,t_2,\ldots,t_{2g+2})})}

\newcommand{\Admuntwo}{\overline{\Adm}({{g_2}\stackrel{2}{\longrightarrow}0})}
\newcommand{\Admunone}{\overline{\Adm}({{g_1}\stackrel{2}{\longrightarrow}0})}
\newcommand{\Admgin}{\overline{\Adm}({g\stackrel{3}{\longrightarrow}\mathbb{P}^1_{\mathbb{C}},((3),t_1,\ldots,t_{2g+2})})}
\newcommand{\Cstar}{\mathbb{C}^\ast}

\newcommand{\Admtinf}{\overline{\Adm}({g\stackrel{3}{\longrightarrow}\mathbb{P}^1_{\mathbb{C}},(t_1,\ldots,t_{2g+4})})}

\newtheorem*{conjecture}{Conjecture}
\newtheorem{lemma}{Lemma}
\theoremstyle{remark}
\newtheorem{defi}[lemma]{Definition}
\restylefloat{figure}

\renewcommand{\to}{\rightarrow}
\let\barrsquare\square
\let\square\undefined

\makeop{Adm}
\makeop{ev}
\makeop{Spec}
\newcommand{\topo}{\mathit{top}}

\begin{document}

\begin{htmlabstract}
<p class="noindent">
In this paper we study a natural class of intersection numbers on
moduli spaces of degree d admissible covers from genus g curves to
<b>P</b><sup>1</sup>, using techniques of localization. These intersection
numbers involve tautological &lambda; and &psi; classes, and are in
some sense analogous to Hodge Integrals on moduli spaces of stable curves.
</p>
<p class="noindent">
We compute explicitly these numbers for all genera in degrees 2 and
3 and express the result in generating function form; we provide a
conjecture for the general degree d case.
</p>
\end{htmlabstract}

\begin{abstract} 
In this paper we study a natural class of intersection numbers on
moduli spaces of degree $d$ admissible covers from genus $g$ curves to
$\mathbb{P}^1$, using techniques of localization. These intersection
numbers involve tautological $\lambda$ and $\psi$ classes, and are in
some sense analogous to Hodge Integrals on moduli spaces of stable curves.

We compute explicitly these numbers for all genera in degrees $2$ and
$3$ and express the result in generating function form; we provide a
conjecture for the general degree $d$ case.
\end{abstract}

\maketitle

\section*{Introduction}

Hodge integrals are a class of intersection numbers on moduli spaces of curves involving   the tautological  classes $\lambda_i$, which are the Chern classes of the Hodge bundle $\mathbb{E}$.
In recent years Hodge integrals have shown a great amount of interconnections with Gromov--Witten theory and enumerative geometry. 

The classical Hurwitz numbers, counting the numbers of ramified Covers of
a curve  with an assigned set of ramification data, can be computed via
Hodge integrals. Simple Hurwitz numbers have been discussed by Ekedahl,
Lando, Shapiro and Vainshtein \cite{elsv:ohnahi,elsv:hnaiomsoc} and
by Graber and Vakil \cite{vg:hnavl}; progress towards double Hurwitz
numbers has been made by Goulden, Jackson and Vakil \cite{gjv:ttgodhn}.

Various spectacular computations of Hodge integrals were carried out in
the late nineties by Faber and Pandharipande \cite{fp:hiagwt}. Their
results have been used to determine the multiple cover contributions
in the GW invariants of $\proj$, thus extending the well-known
Aspinwall--Morrison formula in Gromov--Witten Theory.

Hodge integrals are also at the heart of the theory developed by Bryan
and Pandharipande \cite{bp:tlgwtoc}, studying  the local Gromov--Witten
theory of curves.

It is this last theory that brought our attention to a similar
type of integrals. We study  moduli spaces of admissible covers, a
natural compactification of the Hurwitz scheme. It has been shown by
Abramovich, Corti and Vistoli \cite{acv:ac} that these spaces are smooth
Deligne--Mumford stacks. A class of natural intersection numbers on these
spaces, parallel (and we believe related) to the structure coefficients
of the Topological Quantum Field Theory in Bryan--Pandharipande
\cite{bp:tlgwtoc}, are obtained in the following way. Consider the
diagram of stacks
$$\bfig\barrsquare/->`->``/<1200,500>[\mathcal{U}`\proj`
  \overline{\Adm}
  ({g\stackrel{d}{\longrightarrow}\proj,(\mu_1,t_2,\ldots,t_{n})})`;
  f`\pi``]
  \efig$$
where
\begin{itemize}
\item
$\overline{\Adm}
  ({g\stackrel{d}{\longrightarrow}\proj,(\mu_1,t_2,\ldots,t_{n})})$
denotes the space of (connected) genus $g$, degree $d$, admissible covers
with ramification $(\mu_1,t_2,\ldots,t_{n})$, which we will discuss at
length in \fullref{admcov}.
\item we consider covers that have one arbitrary ramification point
$\mu_1$; all other ramification is simple ($t$ stands for transposition);
\item $\mathcal{U}$ is the universal family;
\item  $f$ is ``morally'' the universal cover map \eqref{f}; 
\end{itemize}

Now define the class of integrals
$$I_d^{\mu}(g) :=
  \int_{\overline{\Adm}
  ({g\stackrel{d}{\longrightarrow}\proj,(\mu_1,t_2,\ldots,t_{n})})}
  \ev_1^\ast(\infty) \cap c_{2g+d-1}(\E),$$
where $\ev_1$ is evaluation at the first marked point \eqref{ev}.

It is an elementary dimension count to show that the only non vanishing
integrals must have $\mu=(d)$, that is, full ramification over the point
$\infty$. For this reason we drop the superscript $\mu$.

We want to organize all of these integrals in generating function form:
$$\mathcal{I}_d(x):= \sum_{g=0}^\infty \frac{I_d(g)}{2g +d-1!}x^{2g+d-1}.$$
These integrals can be approached with techniques of
localization. We follow the spirit, and also the notation, of Faber
and Pandharipande \cite{fp:hiagwt}, who pioneered and developed the
fundamental ideas of auxiliary localization integrals, and of using
different linearizations of line bundles as means to find relations
between Hodge integrals.

We make the following conjecture.

\begin{conjecture}
For all $d \geq 1 $
$$\mathcal{I}_d(x)= (-1)^{d-1}\frac{1}{d}
  \frac{\bigl(2\sin\bigl(\frac{x}{2}\bigr)\bigr)^d}
  {2\sin\bigl(\frac{dx}{2}\bigr)}.$$
\end{conjecture}

The conjecture is trivially true for $d=1$. In this paper we prove it for $d=2,3$. Different strategies are required to prove these two results.  In degree $2$ we exploit the fact that generic ramification is, in fact, full ramification. In degree $3$  we prove the result by means of an auxiliary integral that we know to vanish;  we can obtain the auxiliary integral precisely because full ramification can be  thought of as a degeneration of simple ramification.

The strategy adopted in degree $3$ should in principle work in higher degrees as well. The problem in a direct computation is that
the combinatorial complexity, which is modest in the two cases we examine,  grows dramatically fast.
 
As a corollary of these computations we obtain generating functions for another interesting class of integrals:
$$
J_d(g) := \int_{ \overline{\Adm}({g\stackrel{d}{\longrightarrow}\proj,(t_1,t_2,\ldots,t_{2g+2d-2})})} c_{2g+2d-2}(\F).
$$
In genus $0$,  we recover the Aspinwall--Morrison formula.

\subsection*{Acknowledgements}
I am grateful first and foremost to my advisor, Aaron Bertram, for his
constant support, motivation, and expert guidance. I also thank Y\,P Lee
and Ravi Vakil for carefully listening to my arguments and providing
useful feedback.

\section{Admissible covers}\label{admcov}
Moduli spaces of admissible covers are a ``natural'' compactification of
the Hurwitz scheme. The fundamental idea is that, in order to understand
limit covers, we  allow the base curve to degenerate together with the
cover. Branch points are not allowed to ``come together''; as two or more
branch points tend to collide, a new component of the base curve sprouts
from the point of collision, and  the points transfer onto it. Similarly,
upstairs the cover splits into a nodal cover.

Now more formally: let $(X,p_1,\ldots,p_r)$ be an $r$--pointed nodal curve of genus $g$.
\begin{defi}
An \emph{admissible cover} $\pi\co E\longrightarrow X$ of degree $d$ is a
finite morphism satisfying the following:
\begin{enumerate}
\item $E$ is a connected nodal curve.
\item Every node of $E$ maps to a node of X.
\item The restriction of $\pi\co E{\longrightarrow}X$ to
$X{\setminus}(p_1,\ldots,p_r)$ is \'{e}tale of constant degree $d$.
\item Over a node, locally in analytic coordinates, $X$, $E$ and $\pi$
are described as follows:
\begin{align*}
E\co& e_1e_2=a, \\[-0.5ex]
X\co& x_1x_2=a^n ,\\[-0.5ex]
\pi\co& x_1= e_1^n,\qua  x_2= e_2^n.
\end{align*}
\end{enumerate}
\end{defi}
Moduli spaces of admissible covers were introduced originally by Harris
and Mumford in~\cite{hm:kd}. Intersection theory on these spaces was
for a long time extremely hard and mysterious, mostly because they are
in general not normal, even if the normalization is always smooth. Only
recently in~\cite{acv:ac}, Abramovich, Corti and Vistoli exhibit this
normalization as the stack of balanced stable maps of degree 0 from
twisted curves to the classifying stack $\mathcal{B}S_d$. This way they
attain both the smoothness of the stack and a nice moduli-theoretic
interpretation of it.

We will abuse notation and refer to the Abramovich--Corti--Vistoli spaces
as admissible covers. We will be interested in admissible covers of
$\proj$. In order to estabilish  notation, let us recall our basic
definitions:
\begin{defi} Fix $d\geq1$, and let $\mu_1,\ldots,\mu_n$ be partitions of
$d$. We denote by
$$\overline{\Adm}({g\stackrel{d}{\longrightarrow}0,(\mu_1,\ldots,\mu_n)})$$
the connected component of the stack of  balanced stable maps of degree 0
from a genus 0, $n$--pointed  twisted curve to $\mathcal{B}S_d$ characterized
by the following conditions:
\begin{enumerate}
\item the associated admissible cover (according to the construction
in \cite[page 3566]{acv:ac}) is a nodal curve of genus $g$.
\item let  $x_1,\ldots,x_n$  be the marks on the base curve; the
ramification profile over $x_i$ is required to be of type $\mu_i$.
\end{enumerate}
We call this the stack of admissible covers of degree $d$ and genus $g$
of a genus 0 curve.
\end{defi}

This is either empty or a smooth stack of dimension $n-3=2g+2d+n+\sum
\ell(\mu_i)-nd-5$, where $\ell(\mu_i)$ denotes the length of the partition
$\mu_i$. It admits two natural maps into moduli spaces of curves, as
represented in the following diagram:
$$\bfig
\barrsquare/->`->``/<1200,400>[
  \overline{\Adm}({g\stackrel{d}{\longrightarrow}0,(\mu_1,\ldots,\mu_{n})})`
  \wwbar{M}_g`\wwbar{M}_{0,n}`;```]
  \efig$$ 
In particular, the vertical map has finite fibers.

We also are interested in fixing a parametrization of the base
$\proj$. The objects we parametrize are the same as above,
but the equivalence relation is stricter: we consider two covers
$E_1\rightarrow\proj$, $E_2\rightarrow\proj$ equivalent if there is an
isomorphism $\varphi\co E_1\rightarrow E_2$ that makes the natural triangle
commute. In other words, we are not allowed to act on the base with an
automorphism of $\proj$.

\begin{defi}
We denote by
$$\overline{\Adm}({g\stackrel{d}
  {\longrightarrow}\proj,(\mu_1,\ldots,\mu_n)})$$
the stack of admissible covers of degree d of (a parametrized) $\proj$
by curves of genus g, with n specified branch points having ramification
profile $\mu_1,\ldots,\mu_n$.
\end{defi}
We construct the space of parametrized admissible covers as the stack
of balanced stable maps of degree $d!$ from the category of genus 0,
$n$--pointed twisted curves to the stack quotient $[\proj/S_d]$, where $
S_d$ acts trivially on $\proj$. This is but a slight variation to the
Abramovich--Corti--Vistoli construction. Let us illustrate what happens
over a geometric point $\Spec(\mathbb{C})$:
\begin{figure}[ht!]
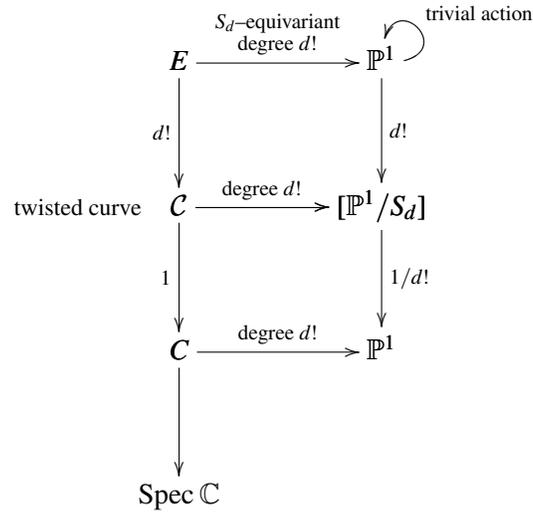

\centering
$$\bfig
  \barrsquare(0,1000)|alra|<700,500>[E`\proj`\mathcal{C}`{[\proj/S_d]};
  \genfrac{}{}{0pt}{1}{S_d\text{--equivariant}}{\text{degree } d!}`
  d!`d!`\text{degree } d!]
  \barrsquare(0,500)|alra|<700,500>[\phantom{\mathcal{C}}`
  \phantom{[\proj/S_d]}`C`\proj;
  `1`1/d!`\text{degree } d!]
  \morphism(0,500)<0,-500>[\phantom{C}`\Spec\mathbb{C};]
  \place(-350,1000)[\text{\footnotesize twisted curve}]
  \morphism(700,1500)|a|/@(r,u)@{->}/<0,0>[\phantom{\proj}`\phantom{\proj};
  \text{trivial action}]
  \efig$$
\caption{The stack of admissible covers of a parametrized $\proj$}
\label{quot1}
\end{figure}

A map of degree $d!$ from the twisted curve produces a map of degree
$1$ from the coarse curve (and this is our desired parametrization
of one special genus 0 twig on the base), a principal  $S_d$ bundle
over the twisted curve  and an $S_d$ equivariant map to $\proj$ (this
data characterizes the admissible cover). Two admissible covers are
equivalent if there is an automorphism of the twisted curve that makes
them commmute. In doing so, the degree $1$ map to $\proj$ has to be
respected, so only the non-parametrized twigs are free to be acted upon
by automorphisms.
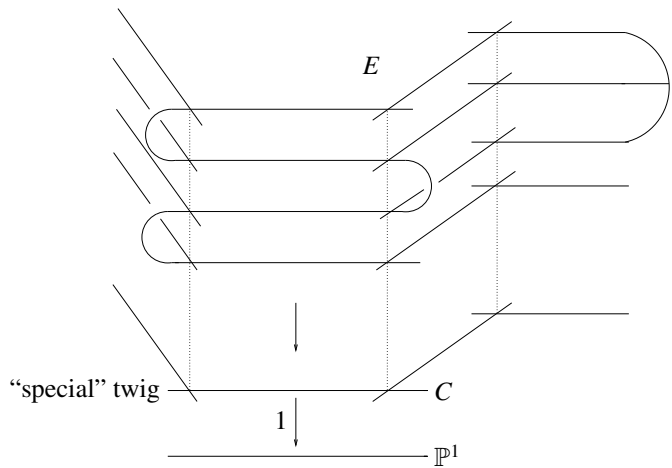
\begin{figure}[ht!]
\centering
\labellist\small
\pinlabel {$E$} at 372 488
\pinlabel {$C$} [l] at 440 87
\pinlabel {$\proj$} [l] at 440 8
\pinlabel {$1$} [r] at 279 52
\pinlabel {``special'' twig} [r] at 122 88
\endlabellist
\includegraphics[width=8cm]{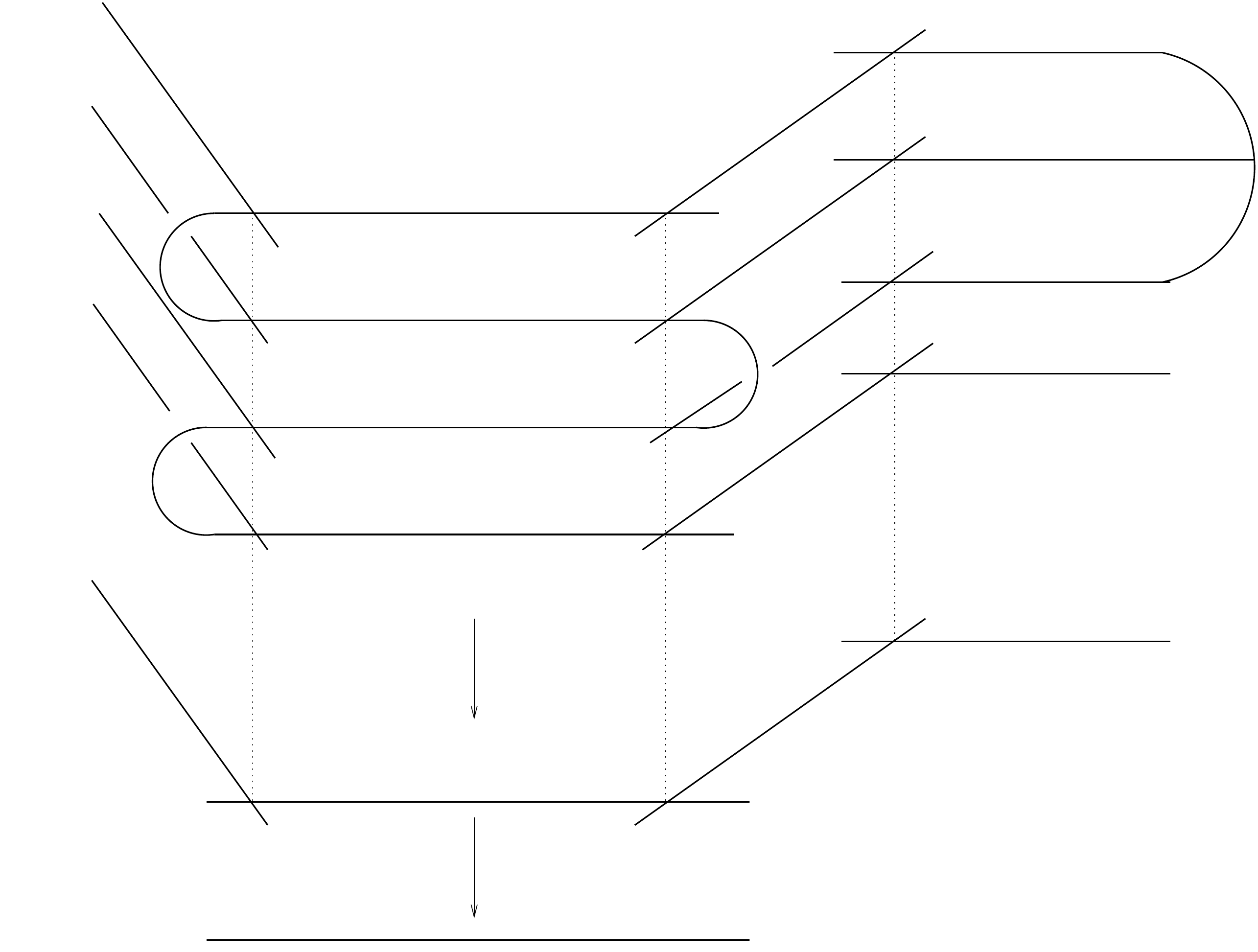}
\label{ciccio}
\caption{Schematic depiction  of an admissible cover of a parametrized $\proj$}
\end{figure} 
This is either empty or a smooth stack of dimension $n=2g+2d+n+\sum
\ell(\mu_i)-nd-2$, admitting two natural morphisms
$$\bfig
  \barrsquare/->`->``/<1200,400>[
  \overline{\Adm}
  ({g\stackrel{d}{\longrightarrow}\proj,(\mu_1,\ldots,\mu_{n})})`
  \wwbar{M}_g`\proj{[n]}`;```]\efig$$
The map to $\wwbar{M}_{g}$ just looks at the source curve forgetting
the cover map. The vertical morphism, taking values in the
Fulton--MacPherson configuration space of $n$ points in $\proj$, looks
instead at the target curve, and at the (ordered) branch points.

The stacks of admissible covers of $\proj$ admit a universal family
$\mathcal{U}$, and a universal cover map $\mu$. The cover map
takes values in a stack $\mathcal{X}$, that is a family over the moduli
space. The fiber over a moduli point consists of a nodal, genus $0$
curve, with one special irreducible component. The universal cover map
can be followed by a  map $\varepsilon$, that contracts  all secondary
twigs and  takes values in $\overline{\Adm} \times\proj$. Finally the
right projection lands us in $\proj$.
\begin{equation}
\bfig
  \ptriangle<500,500>[\mathcal{U}`\mathcal{X}`
  \overline{\Adm}
  ({g\stackrel{d}{\longrightarrow}\proj,(\mu_1,\ldots,\mu_n)});
  \mu`\pi`]
  \morphism(500,500)<1000,0>[\phantom{\mathcal{X}}`
  \overline{\Adm}
  ({g\stackrel{d}{\longrightarrow}\proj,(\mu_1,\ldots,\mu_n)});]\efig
\label{f}
\end{equation}
We call $f$ the composition of the three horizontal maps.

The universal family can be itself interpreted as a moduli space of
admissible covers. If we think of admissible covers as of  stable maps
from a twisted curve, then we obtain a Universal family by adding a mark
to the  twisted curve and requiring trivial ramification over it. Let us
denote with $(1)$ the partition $(1,\ldots,1)$ of $d$, representing an
unramified point. Then,
$$\mathcal{U}= \overline{\Adm}
  ({g\stackrel{d}{\longrightarrow}\proj,(\mu_1,\ldots,\mu_n,(1))}).$$
We can define $n$ tautological sections
$$\sigma_i\co\overline{\Adm}
  ({g\stackrel{d}{\longrightarrow}\proj,
  (\mu_1,\ldots,\mu_n)})\longrightarrow
  \overline{\Adm}({g\stackrel{d}{\longrightarrow}\proj,
  (\mu_1,\ldots,\mu_n,(1))})$$
of the natural forgetful map.
The image of the $i$th section consists of covers where a new rational
component has sprouted from the $i$th marked point. The marked points
$(1)$ and $\mu_i$ have transferred onto this twig. Over this twig we
find $\ell(\mu_i)$ copies of $\proj$ fully ramified over the attaching
point and over the marked point $\mu_i$.

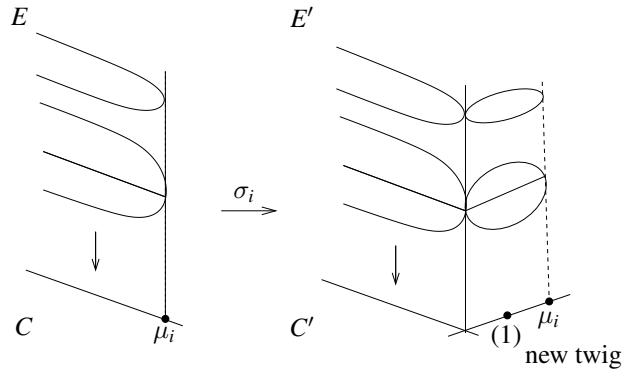
\begin{figure}[ht!]
\centering
\labellist\small
\pinlabel {$E$} at 15 450
\pinlabel {$E'$} at 380 450
\pinlabel {$C$} at 20 50
\pinlabel {$C'$} at 380 50
\pinlabel {$\mu_i$} [t] at 202 60
\pinlabel {$\mu_i$} [t] at 698 80
\pinlabel {$(1)$} [t] at 643 62
\pinlabel {$\sigma_i$} [b] at 305 200
\pinlabel {new twig} at 730 10
\endlabellist
\includegraphics[width=8cm]{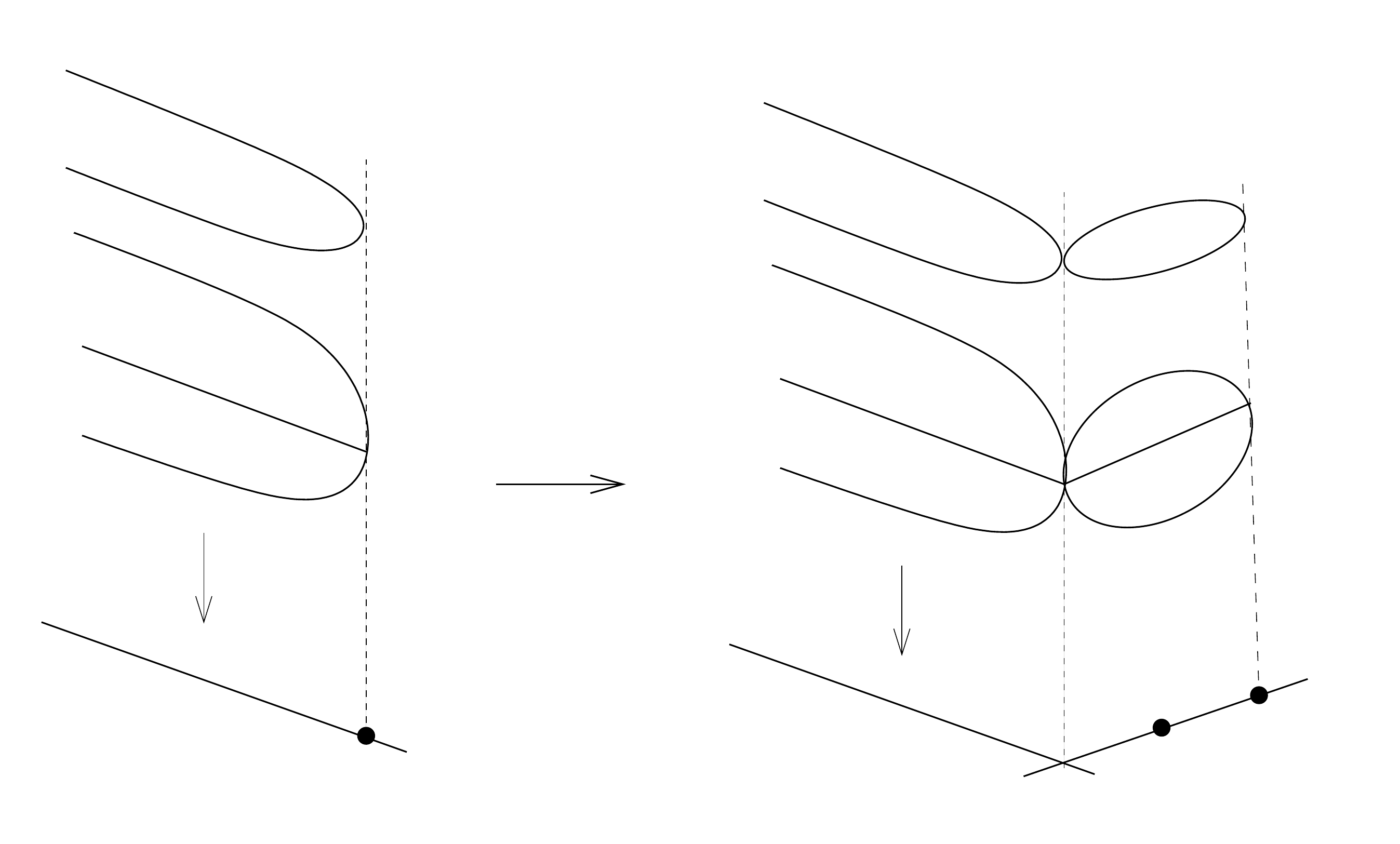}
\label{sigmai}
\caption{The tautological section $\sigma_i$}
\end{figure}

Finally we can define the natural evaluation maps
\begin{equation}
\ev_i := f \circ \sigma_i\co \overline{\Adm}
  ({g\stackrel{d}{\longrightarrow}\proj,(\mu_1,\ldots,\mu_n)})
  \longrightarrow \proj.\label{ev}
\end{equation}

\subsection{The boundary}

The boundary of spaces of admissible covers can be described in terms
of admissible cover spaces of possibly lower degree or genus. In the
case of admissible covers of a parametrized $\proj$, the boundary will
involve also admissible covers of an unparametrized genus $0$ curve. In
\fullref{ciccio}, for example, we can obtain the depicted admissible
cover by ``gluing together'' one admissible cover of a parametrized
$\proj$ (the cover of the special twig) and three admissible covers of
an irreducible genus $0$ curve.  It would be very tempting to conclude
that the irreducible boundary components of an admissible cover space are
actually products of other admissible cover spaces; however we need to
be very careful, and consider the contribution to the stack structure
given by automorphisms.

To illustrate this point let us carefully analyze the gluing map. For
simplicity of exposition, let's  glue at a fully ramified point:
$$\bfig
\barrsquare(-500,0)/`->``^{ (}->/<1400,400>[\phantom{\overline{\Adm}
  (g_1{\stackrel{d}{\longrightarrow}}\proj,
  (\mu_1,\ldots,\mu_{n_1},(d)))}``\mathcal{B}`
  \overline{\Adm}(g_1{+}g_2{\stackrel{d}{\longrightarrow}}\proj,
  (\mu_1,\ldots,\mu_{n_1},\lambda_1,\ldots,\lambda_{n_2}));```]
\morphism(0,400)/{}/<0,0>[\overline{\Adm}
  (g_1{\stackrel{d}{\longrightarrow}}\proj,
  (\mu_1,\ldots,\mu_{n_1},(d))) \times
  \overline{\Adm}(g_2{\stackrel{d}{\longrightarrow}}0,
  (\lambda_1,\ldots,\lambda_{n_2},(d)))`;]
  \efig$$
We claim that the vertical map is an \'{e}tale map of stacks of degree $1/d$. Let us look at a point $[E\rightarrow X]$ of $\mathcal{B}$: we observe that it admits a unique preimage $([E_1\rightarrow \proj],[E_2\rightarrow X_2])$, and we  count the automorphisms of the preimage modulo automorphisms pulled-back from below. In local analytic coordinates around the node, the cover is described as
$$\bfig\morphism<0,-400>[
  \Spec(\mathbb{C}[e_1,e_2]/(e_1e_2-a))`
  \Spec(\mathbb{C}[x_1,x_2]/(x_1x_2-a^d));]\efig$$
by the local equations $x_1=e_1^d, x_2=e_2^d$. Modding out by automorphisms of the ``glued'' cover is equivalent to requiring the first coordinate  $e_1$ to remain untouched. It is then evident that what we have left are $d$ distinct automorphisms, consisting in multiplying $e_2$ by a $d^{th}$ root of unity. This establishes our claim.

Now if we want to glue two branch points with ramification profile
$\eta=(d_1,\ldots,d_{k_\mu})$, with all the $d_i$'s distinct, the situation will be analogous. The gluing map 
\begin{equation}
\bfig
\barrsquare(-800,0)/`->``^{ (}->/<1400,400>[\phantom{X}``\mathcal{B}`
  \overline{\Adm}(g_1{+}g_2{+}\ell(\eta){-}1
  {\stackrel{d}{\longrightarrow}}\proj,
  (\mu_1,\ldots,\mu_{n_1},\lambda_1,\ldots,\lambda_{n_2}));```]
\morphism(0,400)/{}/<0,0>[\overline{\Adm}
  (g_1{\stackrel{d}{\longrightarrow}}\proj,
  (\mu_1,\ldots,\mu_{n_1},\eta)) \times
  \overline{\Adm}(g_2{\stackrel{d}{\longrightarrow}}0,
  (\lambda_1,\ldots,\lambda_{n_2},\eta))`;]
  \efig
\label{glue}
\end{equation}
is an \'{e}tale map of stacks of degree $1/ (d_1\cdot...\cdot d_{k_\mu})$.
 
For the purposes of this paper, this is all we are concerned with. For the sake of a more complete exposition, we briefly describe what the situation for a general partition $\eta$  is. Let
$$\eta=((\eta^1)^{m_1},\dots,(\eta^k)^{m_k}).$$
The gluing map \eqref{glue} is an \'{e}tale map of stacks of degree
\begin{equation}
\frac{1}{\prod (\eta^i)^{m_i}m_i!}.
\label{ionel}
\end{equation}
There is a little bit of combinatorial  subtlety to be dealt with to
obtain this last result. In order to be able to even define a gluing
map we must introduce markings on the covers. Ionel develops the theory
of these spaces in \cite{i:trrih}; our spaces are \'{e}tale quotients
of Ionel spaces. The gluing map is well defined on the level of Ionel
spaces, and it descends to \eqref{ionel}.

\subsection{Tautological classes}

We are interested in describing some ``tautological'' intersection classes
on the stack of admissible covers of an unparametrized genus $0$ curve:
in particular we want to endow our space with analogues  of $\lambda$
and $\psi$ classes. To do so, we will simply pull-back these classes
from the appropriate moduli spaces.

Recall the forgetful map
$$\overline{\Adm}({g\stackrel{d}{\longrightarrow}0, (\mu_1,\ldots,\mu_{n})})
  \stackrel{s}{\longrightarrow}  \wwbar{M}_{g}.$$
The tautological class $\lambda_i \in A^i(\wwbar{M}_{g})$ is defined
to be the $i$th Chern class of the Hodge bundle $\mathbb{E}$.
\begin{defi}
The tautological class $\lambda_i^{\Adm} \in
A^i(\overline{\Adm}({g{\stackrel{d}{\longrightarrow}}0,
(\mu_1,\ldots,\mu_{n})}))$  is defined to be the $i$th Chern class of
the pull-back of the Hodge bundle via the map $s$:
$$\lambda_i^{\Adm}:=s^\ast (\lambda_i).$$
\end{defi}
We will drop the superscript ``$\Adm$'' and simply write $\lambda_i$
whenever there is no risk of confusion.

Let us now look at another natural map
$$\overline{\Adm}({g\stackrel{d}{\longrightarrow}0,
  (\mu_1,\ldots,\mu_{n})}) \stackrel{t}{\longrightarrow}
  \wwbar{\mathcal{M}}_{0,n}.$$
The stack $\wwbar{\mathcal{M}}_{0,n}$ is the moduli space of twisted
$n$--pointed curves of genus 0.

Let $\smash{\wwbar{\mathcal{M}}_{0,n+1}
  \stackrel{\pi}{\longrightarrow}\wwbar{\mathcal{M}}_{0,n}}$
be the universal family over this stack,
{\Large$\omega$}$_{\pi}\longrightarrow\wwbar{\mathcal{M}}_{0,n+1}$ be
the relative dualizing sheaf  and  $\sigma_i$ the $i$th tautological
section. Then $\psi_i \in A^1({\wwbar{\mathcal{M}}_{0,n}})$
is defined to be the first Chern class of $\sigma_i^\ast(\omega_{\pi})$.
\begin{defi}
The  tautological class $\psi_i^{\Adm} \in
A^1(\overline{\Adm}({g\stackrel{d}{\longrightarrow}0,
(\mu_1,\ldots,\mu_{n})}))$  is defined to be the pull-back of the
analogous class via the map $t$:
$$\psi_i^{\Adm}:=t^\ast(\psi_i).$$
\end{defi}
Again, the superscript will be  dropped unless needed for clarity.\\
We can also view $\psi$ classes in a more intrinsic fashion.  Consider
\begin{itemize}
\item the space
$$\overline{\Adm}({g\stackrel{d}{\longrightarrow}0,
  (\mu_1,\ldots,\mu_{n},(1))})$$
where we have  added a trivial ramification condition;
\item the forgetful map
$$\pi_{(1)}\co\overline{\Adm}({g\stackrel{d}{\longrightarrow}0,
(\mu_1,\ldots,\mu_{n},(1))})\longrightarrow
\overline{\Adm}({g\stackrel{d}{\longrightarrow}0, (\mu_1,\ldots,\mu_{n})});$$
\item the $i$th tautological section
$$\sigma_i\co\overline{\Adm}({g\stackrel{d}{\longrightarrow}0,
(\mu_1,\ldots,\mu_{n})})\longrightarrow
\overline{\Adm}({g\stackrel{d}{\longrightarrow}0,
(\mu_1,\ldots,\mu_{n},(1))}).$$
\end{itemize}
\begin{lemma}
The class $-\psi_i\in A^1( \overline{\Adm}({g\stackrel{d}{\longrightarrow}0, (\mu_1,\ldots,\mu_{n})}))$ is the first Chern class of the normal bundle to the image of the section $\sigma_i$.
\end{lemma}
\begin{proof}
Observe the following commutative diagram:
$$\bfig
  \barrsquare|allb|/->`@{->}@<-3pt>`@{->}@<-3pt>`->/<1500,500>[
  \overline{\Adm}
  ({g\stackrel{d}{\longrightarrow}0,(\mu_1,\ldots,\mu_{n},(1))})`
  \wwbar{\mathcal{M}}_{0,n+1}`
  \overline{\Adm}
  ({g\stackrel{d}{\longrightarrow}0,(\mu_1,\ldots,\mu_{n})})`
  \wwbar{\mathcal{M}}_{0,n};
  \tilde{t}`\pi_{(1)}`\pi`t]
  \morphism(0,0)|r|/@{->}@<-3pt>/<0,500>[
  \phantom{\overline{\Adm}
  ({g\stackrel{d}{\longrightarrow}0,(\mu_1,\ldots,\mu_{n},(1))})}`
  \phantom{\overline{\Adm}
  ({g\stackrel{d}{\longrightarrow}0,(\mu_1,\ldots,\mu_{n})})};
  \tilde{\sigma}_i]
  \morphism(1500,0)|r|/@{->}@<-3pt>/<0,500>[
  \phantom{\wwbar{\mathcal{M}}_{0,n+1}}`
  \phantom{\wwbar{\mathcal{M}}_{0,n}};
  \sigma_i]
  \efig$$
We know from Abramovich, Corti and Vistoli \cite[page 3561]{acv:ac}
that the maps $t$ and $\tilde{t}$ are \'{e}tale onto their image. Further,
the diagram is cartesian. Now our lemma follows from the analogous
statement on $\wwbar{\mathcal{M}}_{0,n}$:
$$-\psi_i^{\Adm}= t^\ast(-\psi_i)=c_1( t^\ast\sigma_i^\ast
N_{\sigma_i})=c_1( \tilde{\sigma}_i^\ast \tilde{t}^\ast N_{\sigma_i})=
c_1(\tilde{\sigma}_i^\ast N_{\tilde{\sigma}_i})\proved$$
\end{proof}

\section{Localization}\label{loc}

The main tool for evaluating our integrals is the Atiyah--Bott
localization theorem \cite{ab:tmmaec}.
Consider the  $1$ dimensional algebraic torus $\mathbb{C}^\ast$, and
recall that the $\mathbb{C}^\ast$--equivariant Chow ring of a point is
a polynomial ring in one variable:
$$A^\ast_{\Cstar}(\{\mathit{pt}\},\mathbb{C})= \mathbb{C}[\hbar]$$
Let $\Cstar$ act on a smooth, proper stack $X$, denote by
$i_k\co F_k\hookrightarrow X$ the irreducible components of the fixed locus
for this action and by $N_{F_k}$ their normal bundles. The natural map:
\begin{eqnarray*}
A^\ast_{\Cstar}(X) \otimes \mathbb{C}(\hbar)
  &\longrightarrow&
  \sum_{k}{A^\ast_{\Cstar}}(F_k) \otimes \mathbb{C}(\hbar)\\
\alpha &\longmapsto& \frac{i_k^\ast\alpha}{c_{\topo}(N_{F_k})}.
\end{eqnarray*}
is an isomorphism. Pushing forward equivariantly to the class of a point,
we obtain the Atiyah--Bott integration formula
$$\int_{[X]}\alpha = \sum_k \int_{[F_k]}
  \frac{i_k^\ast\alpha}{c_{\topo}(N_{F_k})}.$$

\subsection{Our set-up}
Let $\mathbb{C}^\ast$ act on a 2--dimensional vector space $V$ via
$$t\cdot(z_0,z_1)=(tz_0,z_1).$$
This action descends on $\proj$, with fixed points $0=(1:0)$ and
$\infty=(0:1)$. An equivariant lifting of $\Cstar$ to a line bundle $L$
over $\proj$ is uniquely determined by its  weights $\{L_0,L_\infty\}$
over the fixed points.

The canonical lifting of $\Cstar$ to the tangent bundle of $\proj$
has weights $\{1,-1\}$.

The action on $\proj$ induces an action on the moduli spaces of admissible
covers to a parametrized $\proj$ simply by postcomposing the cover map
with the automorphism of $\proj$ defined by $t$.

The fixed loci for the induced action on the moduli space consist of
admissible covers such that anything ``interesting'' (ramification,
nodes) happens over $0$ and $\infty$, or on ``non-special'' twigs that
attach to the main $\proj$ at  $0$ or $\infty$.

\subsection{Restricting Chow classes to the fixed loci}
We want to compute the restriction to various fixed loci of the top
Chern class of the bundle
$$E=R^1\pi_\ast f^\ast(\mathcal{O}_{\proj}\oplus\mathcal{O}_{\proj}(-1)).$$
The top Chern class $c_{2g+d-1}(E)$ splits as
$$c_{2g+d-1}(E)=c_g(R^1\pi_\ast f^\ast\mathcal{O}_{\proj})
  c_{g+d-1}(R^1\pi_\ast f^\ast\mathcal{O}_{\proj}(-1)),$$
so we will analyze the two terms separately.

There is a standard technique to carry out these computations. To avoid
an overwhelmingly cumbersome notation, we choose to show it only in a
particular example, that will be the most important for our purposes.

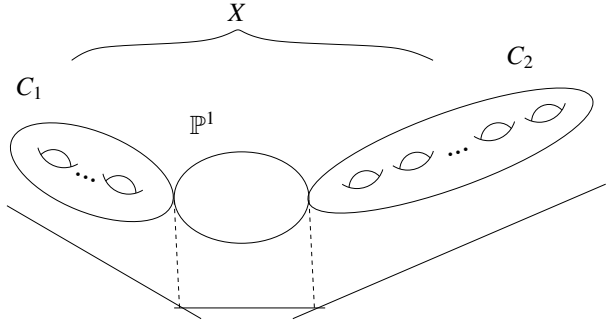
\begin{figure}[ht!]
\labellist\small
\pinlabel {$X$} [b] at 193 243
\pinlabel {$C_1$} at 20 194
\pinlabel {$C_2$} at 430 220
\pinlabel {$\proj$} at 165 165
\endlabellist
\includegraphics[width=8cm]{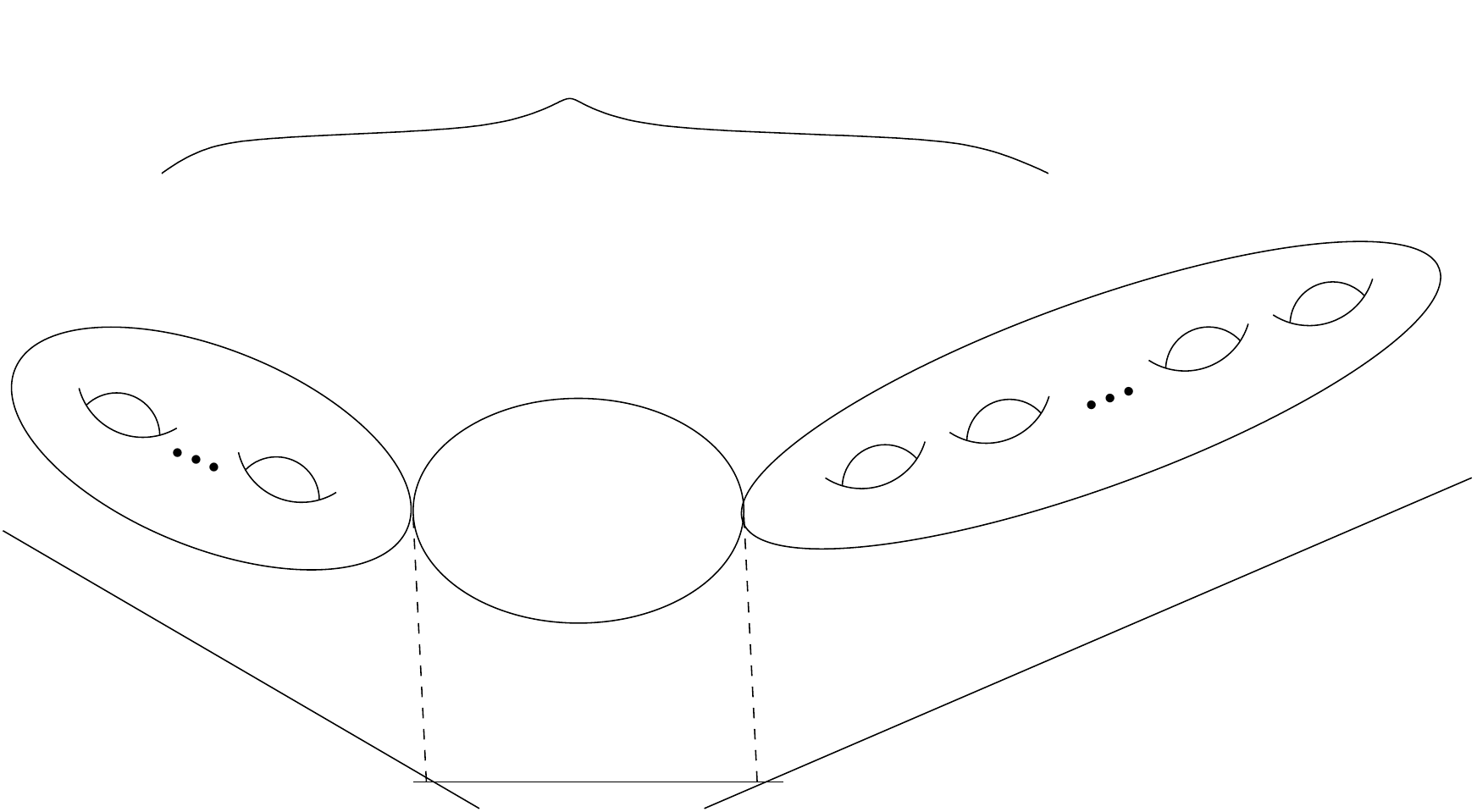}
\caption{The fixed locus $F_{g_1,g_2}$}
\label{f12}
\end{figure}

Let's consider the fixed locus $F_{{g_1}{g_2}}$, consisting of covers
where the main $\proj$ is ramified over $0$ and $\infty$ and curves of
genus $g_1$ and $g_2$ are attached on either side. A point in this fixed
locus is represented in \fullref{f12}, where we denote by $X$ the nodal
curve, $C_1$ and $C_2$ the irreducible components over $0$ and $\infty$.

The starting point in analyzing the restriction of the bundle $E$ to this fixed locus  is the classical \textit{normalization sequence}:
$$0\longrightarrow\mathcal{O}_{X}\longrightarrow\mathcal{O}_{C_1}\oplus\mathcal{O}_{\proj}\oplus\mathcal{O}_{C_2}\longrightarrow\mathbb{C}_{n_1}\oplus\mathbb{C}_{n_2}\longrightarrow 0.$$

\textbf{Term 1} ($c_g(R^1\pi_\ast f^\ast\mathcal{O}_{\proj})$)\qua
It suffices to analyze the long exact sequence in cohomology associated
to the normalization sequence
\begin{multline*}
0\longrightarrow h^0(\mathcal{O}_{X})\longrightarrow
  h^0(\mathcal{O}_{C_1})\oplus h^0(\mathcal{O}_{\proj})\oplus
  h^0(\mathcal{O}_{C_2})\longrightarrow
  \mathbb{C}_{n_1}\oplus\mathbb{C}_{n_2} \\
\longrightarrow h^1(\mathcal{O}_{X})\longrightarrow
  h^1(\mathcal{O}_{C_1})\oplus h^1(\mathcal{O}_{C_2})\longrightarrow 0.
\end{multline*}
Assume that $\mathcal{O}_{\proj}$ is linearized with weights $\{\alpha,\alpha\}$.
Then
$$c_g(R^1\pi_\ast f^\ast\mathcal{O}_{\proj})= (-)^g \Lambda_{g_1}(-\alpha) \Lambda_{g_2}(-\alpha),$$
where the following notational convention holds:
$$\Lambda_g(n)=\sum (n\hbar)^i\lambda_{g-i}.$$
The reason for switching from $\alpha$ to $-\alpha$ is that $h^1(\mathcal{O})$ are the fibers of the dual bundle to the Hodge bundle, hence the odd degree Chern classes will have a negative sign.

\textbf{Term 2} ($c_{g+d-1}(R^1\pi_\ast f^\ast\mathcal{O}_{\proj}(-1))$)\qua
In this case we first want to tensor the normalization sequence by $
f^\ast\mathcal{O}_{\proj}(-1)$, and then proceed to analyze the long
exact sequence in cohomology:
\begin{multline*}
0\longrightarrow h^0(\mathcal{O}_{C_1})\oplus h^0(\mathcal{O}_{C_2})
  \longrightarrow\mathbb{C}_{n_1}\oplus\mathbb{C}_{n_2} \\
\longrightarrow h^1( f^\ast\mathcal{O}_{\proj}(-1)) \longrightarrow
  h^1(\mathcal{O}_{C_1})\oplus h^1(\mathcal{O}_{\proj}(-d))\oplus
  h^1(\mathcal{O}_{C_2})\longrightarrow 0.
\end{multline*}
Now, having linearized $\mathcal{O}_{\proj}(-1)$ with weights $\{\beta,\beta+1\}$, 
$$c_{g+d-1}(R^1\pi_\ast
  f^\ast\mathcal{O}_{\proj}(-1))=(-)^g \Lambda_{g_1}(-\beta)
  \Lambda_{g_2}(-\beta-1)\hbar^{d-1}\prod_1^{d-1}
  \Bigl(\beta+\frac{i}{d}\Bigr).$$
The last term in our contribution, coming from
$h^1(\mathcal{O}_{\proj}(-d))$, is explained in the following
way. Consider a degree $d$ map from $\proj$ to $\proj$. The target curve
is given the natural $\mathbb{C}^\ast$ action, and the tautological
bundle is linearized with weights $\{\beta,\beta+1\}$. Now let $x$ and
$z$ be local coordinates around $0$ for, respectively, the target and
the source curve. The expression of the map in local coordinates is
$$x=z^d.$$
We see then that $z$ must have weight $-1/d$. The vector space
$h^1(\mathcal{O}_{\proj}(-d))$ is $(d-1)$ dimensional and generated, in
local coordinates, by the sections
$\{1/z,1/z^2,\ldots,1/z^{d-1}\}$. The
line bundle over moduli with these fibers is trivial, because $\proj$
is rigid, but it is linearized with weights $\beta +i/d$; $\beta$
coming from the weight of the trivialization of the pullback of
$\mathcal{O}_{\proj}(-1)$ in the chart over $0$, $i/d$ from the section
$1/z^i$. Notice that if you were to reproduce this computation using a
local cohordinate over $\infty$ instead, the corresponding weights would
now be $(\beta+1)-(d-i)/d$, which are exactly the same.

\subsection{The Euler class of the normal bundle to the fixed loci}

The standard way to carry out this computation is to analyze the
deformation long exact sequence, and identify the fiber of the normal
bundle to a fixed locus at a particular moduli point to the moving part
(the part where the $\Cstar$ action doesn't lift trivially) of the
tangent space to the moduli space (corresponding to the space of first
order deformations of the admissible cover in question). It's shown
by Abramovich, Corti and Vistoli \cite[page 3561]{acv:ac} that the
deformation theory of admissible covers corresponds exactly to the
deformation theory of the base, genus $0$, twisted curve. The reason for
this is that admissible covers are  \'{e}tale covers (in fact principal
$S_d$--bundles) of the base twisted curve.

Deformations of a genus $0$ nodal twisted curve are described as follows:
first of all, we can deal with one node at a time. For one given node,
there are two different potential contributions:
\begin{itemize}
\item the contribution from  moving the node on the main $\proj$. Doing
this infinitesimally means moving along the tangent space to the
attaching point on the main $\proj$. Again, the bundle with fiber the
tangent space over a given point of $\proj$ is a trivial bundle, but
in equivariant cohomology it can have a purely equivariant first Chern
class, according to the linearization of the fibers. In our particular
case, the tangent bundle has weight $1$ over $0$ and $-1$ over $\infty$,
thus producing a contribution of $\hbar$ for moving a node around $0$,
of $-\hbar$ for moving a node around $\infty$;
\item the contribution from smoothing the node. It corresponds to
the first Chern class of the tensor product of the tangent spaces at
the attaching points  of the two curves. Again, we get a $\pm\hbar$
contribution from the point on the main $\proj$; the other attaching point
$x$, on the other hand, contributes, by definition, a $-\psi_x$ class.
\end{itemize}
  
\section{Degree 2}
We now carry out the explicit computation of the integral
$$I_2(g)=\int_{\Admtwog}\ev_1^\ast(\infty)\cap c_{2g+1}(\E),$$
for all genera, and express the result in generating function form:
$$\mathcal{I}_2(x)= \sum_{g=0}^\infty
  \Bigl(\frac{I_2(g)}{2g+1!}\Bigr)x^{2g+1}.$$

\subsection{The strategy}
It is important to notice that, while the final result is independent
of the choice of the lifting of the $\mathbb{C}^{\ast}$ action to the
vector bundle $E=\E$, the intermediate calculations are not. This is
in fact the heart of our strategy. We choose two different specific
linearizations with the twofold objective of
\begin{itemize}
\item limiting a priori the number and the combinatorial complexity of
the contributing fixed loci;
\item obtaining, by equating the calculations with the two linearizations,
a recursive formula for genus $g$ integrals in term of lower genus data.
\end{itemize}

\subsection{The localization set-up}

We induce different linearizations on the bundle $E$ by choosing different
liftings of the $\Cstar$ action on the bundles $\mathcal{O}_{\proj}$
and $\mathcal{O}_{\proj}(-1)$. Recall that a linearization of a line
bundle over $\proj$ is determined by the weights of the fixed fibers
representations.

\textbf{Linearization A}\qua We choose to linearize the two bundles as
indicated in the following table:
\begin{center}
\begin{tabular}{|l||c|c|}
\hline
weight & over $0$ & over $\infty$ \\
\hline
\hline
 ${\mathcal{O}_{\proj}}(-1)$ & $-1$ & 0 \\
\hline
${\mathcal{O}_{\proj}}$  & 0  & 0  \\
\hline
\end{tabular}
\end{center}
There is only one fixed locus ${F_{g,\cdot}}$ contributing to
the localization integral, consisting in a cover of $\proj$ fully
ramified over $0$ and $\infty$, and a genus $g$ curve mapping with
degree $2$ to an unparametrized $\proj$ sprouting from the point
$0$. \fullref{{F_{g,cdot}}} illustrates the fixed locus, and the
conventional graph notation to indicate it.

\begin{figure}[ht!]
\centering
\labellist\small
\pinlabel {$C_1$} at 40 220
\pinlabel {$\proj$} at 190 190
\pinlabel {genus 0 twig} [tr] at 94 80
\pinlabel {$\proj$} at 290 15
\pinlabel {$g$} [b] at 437 107
\endlabellist
\includegraphics[width=10cm]{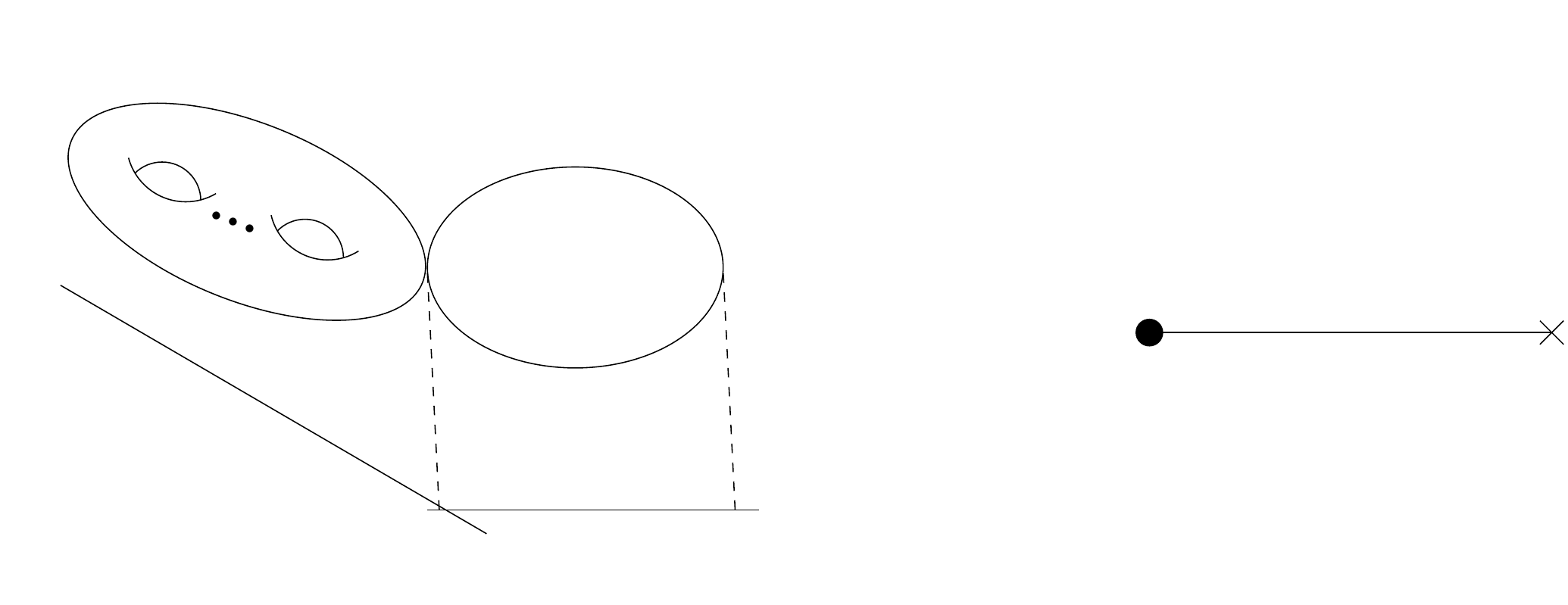}
\caption{The fixed locus ${F_{g,\cdot}}$}
\label{{F_{g,cdot}}}
\end{figure}
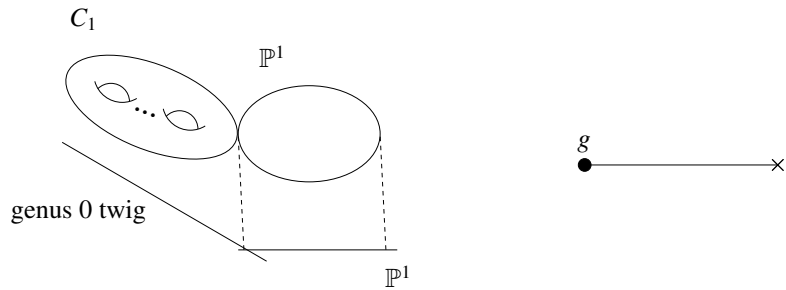

The reason for this dramatic collapsing of the contributing fixed loci
lies in some standard localization facts:
\begin{itemize}\label{fl}
\item the  ramification condition required over $\infty$ implies
that there can be only one connected component in the preimage of
$\infty$. This translates to the fact that the localization graph can
have at most $1$ vertex over $\infty$;
\item the weight $0$ linearization of ${\mathcal{O}_{\proj}}(-1)$ over
$\infty$ implies that the localization graph must have valence $1$
over $\infty$.
\item finally, let's observe that both bundles have weight $0$
over $\infty$; the restriction of our bundle to fixed loci that
have contracted components over $\infty$  involves the class
$\lambda_{g_\infty}^2$, that vanishes  for $g>0$ by a famous result by
Mumford \cite{m:taegotmsoc}. The only option is then to have genus $0$
over infinity. But a genus zero curve with only two special  points is
instable, and hence must be contracted.
\end{itemize}

\textbf{Linearization B}\qua We choose to linearize the bundles  with weights:
\begin{center}
\begin{tabular}{|l||c|c|}
\hline
weight & over $0$ & over $\infty$ \\
\hline
\hline
 ${\mathcal{O}_{\proj}}(-1)$ & $-1$ & 0 \\
\hline
${\mathcal{O}_{\proj}}$  & 1  & 1  \\
\hline
\end{tabular}
\end{center}

In this case the analysis of the possibly contributing fixed loci is similar, except we can't appeal to Mumford's relation any more. Hence our fixed loci will consist of a copy of $\proj$ ramified over $0$ and $\infty$, with two curves of genus $g_1$,$g_2$ attached on either side. (And, of course, $g_1 + g_2 = g$). These are the loci $F_{g_1,g_2}$ described in \fullref{f12}.

\subsection{Explicit evaluation of the integral and recursion}

\textbf{Linearization A}\qua
Let us first of all observe that ${F_{g,\cdot}}$ is naturally isomorphic
to $\smash{\Admtwoun}$.  Using the computations in \fullref{loc},
and the standard equivariant cohomology fact that $\ev_1^\ast(\infty)=
-\hbar$, we obtain the explicit evaluation of our integral on this
fixed locus:
\begin{multline*}
I_2^A(g)=\int_{\Admtwoun}\frac{\lambda_g\Lambda_g(1)
  (-\hbar/2)}{\hbar(\hbar-\psi)} \\
= -\frac{1}{2}\int_{\Admtwoun}
  \lambda_g\lambda_{g-1}+\lambda_g\lambda_{g-2}\psi+
  \cdots+\lambda_g\psi^{g-1}.
\end{multline*}
Just as a convenient notation, let's denote the last integral by $L_2(g)$,
so that
\begin{equation}
I_2^A(g)=-\tfrac{1}{2}L_2(g).
\label{norm}
\end{equation}
\textbf{Linearization B}\qua
In this case we have  $g+1$ different types of fixed loci, corresponding
to all possible ways of choosing an ordered pair of nonnegative integers
adding to $g$. We will study separately three situations:
\begin{description}
\item[$F_{\cdot,g}$] This fixed locus is naturally isomorphic to $2g+1$
disjoint copies of
$$\Admtwoun.$$
The evaluation of the integral reads
\begin{multline*}
\int_{F_{\cdot,g}} \frac{\lambda_g\Lambda_g(-1)
  (-\hbar/2)}{(-\hbar-\psi)}= \\
-\frac{2g+1}{2}\int_{\Admtwoun}\lambda_g\lambda_{g-1}+
  \lambda_g\lambda_{g-2}\psi + \cdots
  +\lambda_g\psi^{g-1} \\
=-\frac{2g+1}{2}L_2(g).
\end{multline*}

\item[$F_{g_1,g_2}$, $g_1,g_2\neq 0)$] After keeping track of the
combinatorics of the gluing  and of the possible distributions of the
marks, the integral evaluates
\begin{align*}
\int_{F_{g_1,g_2}}
&\frac{\Lambda_{g_1}(-1)\Lambda_{g_1}(1)\lambda_{g_2}\Lambda_{g_2}(-1)
(-\hbar/2)}{\hbar(\hbar-\psi)(-\hbar-\psi)} \\
&=-{\binom{2g+1}{2g_2}}\int_{\Admunone}(-)^{g_1}\psi^{2g_1-1} \\
&\hspace{90pt}\int_{\Admuntwo}
  \lambda_{g_2}\lambda_{g_2-1}+ \lambda_{g_2}\lambda_{g_2-2}\psi + \cdots
  +\lambda_{g_2}\psi^{g_2-1} \\
&\hspace{165pt}:=(-)^{g_1+1}\frac{1}{2}{{2g+1}\choose{2g_2}}P_2(g_1)L_2(g_2).
\end{align*}
To make the notation a little lighter we omitted the marked points
(that are still there, though). Also we choose to denote with $P_2$
the integral of $\psi_t$ to the top power.

\item[$F_{g,\cdot}$] This is the same fixed locus encountered in the
computations with linearization A. However, the contribution in this
case will be quite different:
\begin{multline*}
\int_{F_{g,\cdot}}
  \frac{\Lambda_g(1)\Lambda_g(-1)(-\hbar/2)}{(\hbar-\psi)} \\
=-\tfrac{1}{2}\int_{\Admtwoun}(-1)^g\psi^{2g-1}=(-)^{g+1}\tfrac{1}{2}P_2(g).
\end{multline*}
\end{description}

So, altogether, the integral computed with Linearization B is
$$I_2^B(g)=-\tfrac{1}{2}(2g+1)L_2(g)-\sum_{i=0}^{g-1}
  (-1)^{g-i}{\binom{2g+1}{2i}}P_2(g-i)L(i),$$
where we have incorporated the last contribution in the summation by
defining $L(0)=1/2$.

\begin{lemma}
For any $i$, $P(i)=\frac{1}{2}$.
\end{lemma}

\begin{proof}
This follows easily from the fact that the $\psi$ classes that we
are using are pulled back on the space of admissible covers  from
$\wwbar{\mathcal{M}}_{0,2g+2}$ via an \'{e}tale,  degree $1/2$
map (this accounts for the hyperelliptic involution upstairs). The
projective coarse moduli space of $\wwbar{\mathcal{M}}_{0,2g+2}$
is $\wwbar{M}_{0,2g+2}$, and the two spaces are birational. It is
a classical result that the integral of $\psi$ to the top power on
$\wwbar{M}_{0,2g+2}$ is one, hence the lemma.
\end{proof}

We can now equate the results obtained with the two different
linearizations, to obtain a recursive formula for the $L_2(g)$'s.
$$-\tfrac{1}{2}L_2(g)=-\tfrac{1}{2}(2g+1)L_2(g)-\sum_{i=0}^{g-1}
  (-1)^{g-i}{\binom{2g+1}{2i}}P_2(g-i)L_2(i)$$
After a tiny bit of elementary arithmetic we obtain
\begin{equation}
%\begin{imp}
L_2(g)= \frac{1}{2g}\sum_{i=0}^{g-1} (-)^{g-i+1}{\binom{2g+1}{2i}}L_2(i).
\label{rel2}
%\end{imp}
\end{equation}

\subsection{The generating function}
We now want to use relation \eqref{rel2} to compute the  generating
function:
$$\mathcal{L}_2(x)=\sum_{i=0}^{\infty}\left(\frac{{L_2}(i)}{2i+1!}\right)
  x^{2i+1}.$$
Let us first of all differentiate this function,
$$\frac{d}{dx}\mathcal{L}_2(x)=
  \sum_{i=0}^{\infty}\left(\frac{{L_2}(i)}{2i!}\right) x^{2i}.$$
Now let us compute
\begin{align*}
\frac{d}{dx}\mathcal{L}_2(x)\cdot{\sin{(x)}}
&=\sum_{g=0}^{\infty}x^{2g+1}\sum_{i=0}^{g}(-)^{g-i+1}
  \frac{{L_2}(i)}{2i!(2g-2i+1)!} \\
&=\sum_{g=0}^{\infty}x^{2g+1}\sum_{i=0}^{g}(-)^{g-i+1} {\binom{2g+1}{2i}}
  \frac{{L_2}(i)}{2g+1!} \\
&=\sum_{g=0}^{\infty}x^{2g+1}\left(\frac{{L_2}(g)}{2g+1!}\right) \\
&=\mathcal{L}_2(x).
\end{align*}
Hence  relation \eqref{rel2} translates to the following ODE on the generating function $\mathcal{L}_2(x)$:
\begin{eqnarray}
%\begin{imp}
{\mathcal{L}}^\prime_2(x)\cdot \sin{(x)}=\mathcal{L}_2(x), \qquad
{\mathcal{L}}_2(0)=0.
%\end{imp}
\label{explicitode}
\end{eqnarray}

This equation integrates to give us $\mathcal{L}_2(x)=
\tan(x/2)$. Finally, recalling \eqref{norm} we can conclude
\begin{eqnarray}
%\begin{vimp}
\mathcal{I}_2(x)=- \tfrac{1}{2}\tan\left(\tfrac{x}{2}\right).
\label{gf2}
%\end{vimp}
\end{eqnarray}

\subsection{A corollary}
Using result \eqref{gf2} it's now easy to compute the generating function for the second class of integrals we are interested in. Consider
$$J_2(g)=\int_{\Admtwog} c_{2g+2}(\F),$$
and the corresponding generating function
$$\mathcal{J}_2(x)= \sum_{g=0}^\infty
  \left(\frac{J_2(g)}{2g+2!}\right)x^{2g+2}.$$
Again, there is a particularly favorable choice of linearizations:

\begin{center}
\begin{tabular}{|l||c|c|}
\hline
weight & over $0$ & over $\infty$ \\
\hline
\hline
 ${\mathcal{O}_{\proj}}(-1)$ & $-1$ & 0 \\
\hline
${\mathcal{O}_{\proj}}(-1)$  & 0  & 1  \\
\hline
\end{tabular}
\end{center}
The only contributing fixed loci must have valence $1$ both over $0$ and
$\infty$. These are precisely the loci $F_{g_1,g_2}$ studied above. The
explicit computation of the integral is
\begin{align*}
\labellist\small \pinlabel {$g$} [b] at 31 25 \endlabellist
\includegraphics[width=3cm]{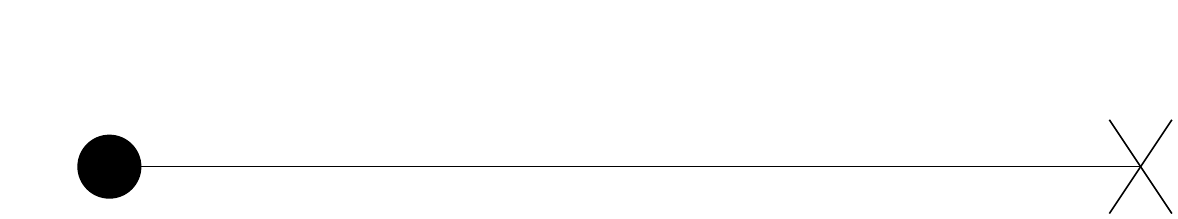} &:
  (2g+2)\int\frac{\lambda_g\Lambda_g(1)\bigl(\frac{\hbar}{2}\bigr)
  \bigl(-\frac{\hbar}{2}\bigr)}{\hbar(\hbar-\psi)(-\hbar)}
= \tfrac{1}{4}(2g+2)L_2(g)\\
\labellist\small
\pinlabel {$g_1$} [b] at 40 20
\pinlabel {$g_2$} [b] at 338 20
\endlabellist
\includegraphics[width=3cm]{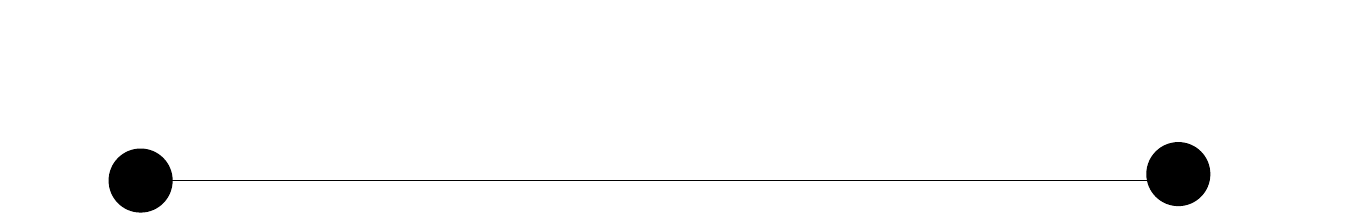} &:
  2{\binom{2g+2}{2g_1+1}}\int\frac{\lambda_{g_1}\Lambda_{g_1}(1)
  \bigl(\frac{\hbar}{2}\bigr)}{\hbar(\hbar-\psi)}
  \int\frac{\lambda_{g_2}\Lambda_{g_2}(-1)
  \bigl(-\frac{\hbar}{2}\bigr)}{-\hbar(-\hbar-\psi)} \\
&\hspace{150pt}= \frac{1}{2}{{2g+2}\choose{2g_1+1}}L_2(g_1)L_2(g_2)\\
\labellist\small \pinlabel {$g$} [b] at 308 26 \endlabellist
\includegraphics[width=3cm]{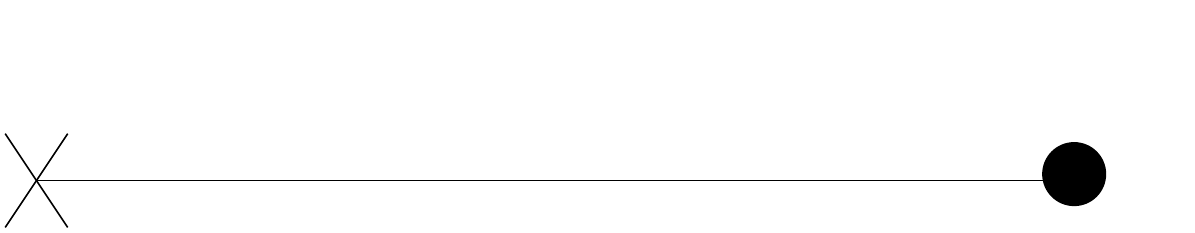} &:
  (2g+2)\int\frac{\lambda_g\Lambda_g(-1)
  \bigl(\frac{\hbar}{2}\bigr)\bigl(-\frac{\hbar}{2}\bigr)}
  {\hbar(-\hbar-\psi)(-\hbar)}
= \tfrac{1}{4}(2g+2)L_2(g)
\end{align*}
All previous integrals are computed over the appropriate unparmetrized
admissible cover spaces. Adding everything together we obtain the
relation
\begin{equation}
J(g)= \frac{1}{2}\sum_0^g{{2g+2}\choose{2i+1}}L_2(i)L_2(g-i).
\end{equation}
(Recalling that we have defined $L_2(0)=1/2$.)

This relation allows us to obtain the generating function
$\mathcal{J}_2(x)$. For this purpose it suffices to notice
\begin{equation}
%\begin{vimp}
\mathcal{J}_2(x)= 2 \mathcal{I}_2(x)^2=
  \tfrac{1}{2}\tan^2\bigl(\tfrac{x}{2}\bigr).
%\end{vimp}
\end{equation}

\section{Degree 3}
In this section we will compute the integral
$${I}_3(g):=\int_{\Admgin}
  \ev_{(3)}^\ast(\infty)\cap c_{2g+2}(\E),$$
for all genera $g$, and present the result in generating function form
$${\mathcal{I}}_3(x):=\sum_{g=0}^\infty \frac{I_3(g)}{2g+2!} x^{2g+2}.$$

\subsection{The strategy}
We will use localization to compute our integral. First of all, we choose
an extremely convenient choice of linearizations on the $\proj$--bundles
${\mathcal{O}_{\proj}}$ and ${\mathcal{O}_{\proj}}(-1)$. This will express
our integral in terms of a Hodge integral over only one boundary component
of the moduli space.

We then will introduce an auxiliary integral, that we know to vanish
for elementary dimension considerations. Evaluating this integral via
localization will produce relations between the integrals $I_3(g)$, for
different genera $g$,  integrals in degree $2$ and simple Hurwitz numbers.

We are able to transform these relations into a linear differential
equation for the generating function ${\mathcal{I}_3}(x)$. Finally,
solving the ODE with the appropriate boundary conditions gives us
the result.

\subsection{The localization set-up}

We choose to linearize the bundle as in linearization A in the previous
section:

\begin{center}
\begin{tabular}{|l||c|c|}
\hline
weight & over $0$ & over $\infty$ \\
\hline
\hline
 ${\mathcal{O}_{\proj}}(-1)$ & $-1$ & 0 \\
\hline
${\mathcal{O}_{\proj}}$  & 0  & 0  \\
\hline
\end{tabular}
\end{center}

For completely analogous reasons to the degree two case (see
\fullref{fl}), there is only one fixed locus, $F_{g,\cdot}$, contributing
to the localization integral, consisting in a cover of $\proj$ fully
ramified over $0$ and $\infty$, and a genus $g$ curve mapping with degree
$3$ to an unparametrized $\proj$ sprouting from the point $0$.

The integral then becomes
\begin{align*}
I_3(g) &= \int_{F_{g,\cdot}}\frac{\lambda_g\Lambda_g(1)
\frac{2}{9}\hbar^2}{\hbar(\hbar-\psi_3)} \\
&= \frac{2}{9}\int_{\Admun} \lambda_g\lambda_{g-1}\psi_3 +
\lambda_g\lambda_{g-2}\psi_3^2+\cdots+ \lambda_g\psi_3^{g}.
\end{align*}
With the sole purpose of keeping track of  coefficients in a more natural
way in what follows, we give a name to the rightmost integral without
the $\frac29$ in front of it:
$$L_3(g):=\int_{\Admun}
  \lambda_g\lambda_{g-1}\psi_3 + \lambda_g\lambda_{g-2}\psi_3^2+\cdots+
  \lambda_g\psi_3^{g}.$$

\subsection{The auxiliary integral}
Let us now consider the following equivariant integral:
$$\int_{\Admtinf}\ev_1^\ast(\infty)\cap c_{2g+2}(\E).$$
This integral must vanish for dimension reasons. Let us now evaluate
this integral via localization.  We now choose different linearizations
for the two bundles, as indicated in the following table.
\begin{center}
\begin{tabular}{|l||c|c|}
\hline
weight & over $0$ & over $\infty$ \\
\hline
\hline
 ${\mathcal{O}_{\proj}}(-1)$ & $-1$ & 0 \\
\hline
${\mathcal{O}_{\proj}}$  & 1  & 1  \\
\hline
\end{tabular}
\end{center}
With this choice of linearizations, the explicit evaluation of the
integral follows. We will again be invoking a famous relation by Mumford
\cite{m:taegotmsoc}:
$$\Lambda_g(-1)\Lambda_g(1)= (-)^g \hbar^{2g}.$$ 

\begin{description}
\item[$F_{g,0}$]
\labellist\small
\pinlabel {$g$} [b] at 43 26
\pinlabel {$3$} [b] at 223 13
\pinlabel {$0$} [b] at 416 22
\endlabellist
\includegraphics[width=3cm]{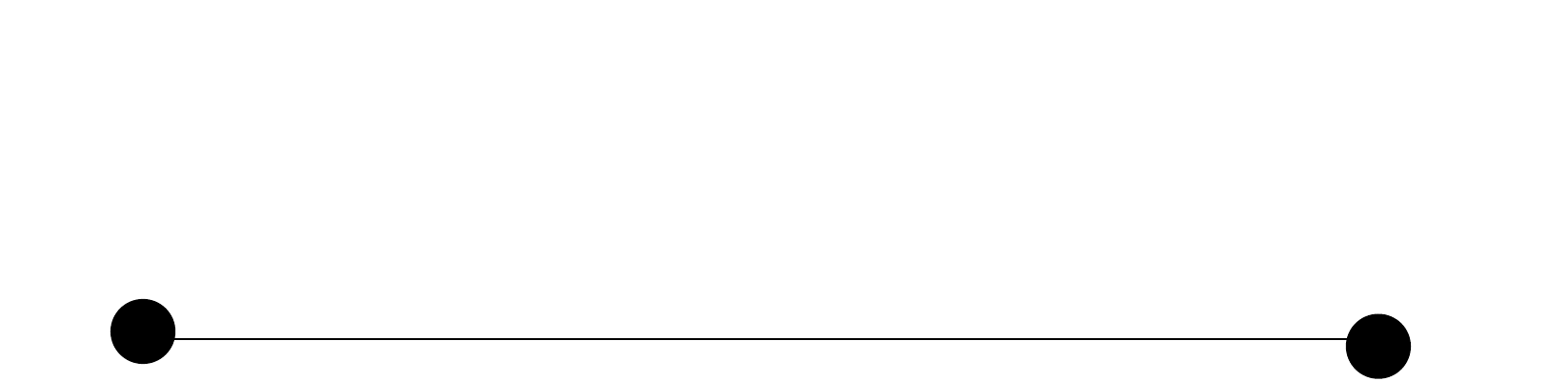}
\begin{multline*}
3{{2g+3}\choose{2g+2}}\int_{\small{{\Admun}}}
{\frac{(-)^g\hbar^{2g}}{\hbar(\hbar-\psi_3)}} \int_{\Admunz}
^{\frac{1}{-\hbar-\psi_3}\left(\frac{2}{9}\hbar^2\right)} \\
= (-)^{g+1}\frac{2}{3}{{2g+3}\choose{2g+2}}
  \frac{1}{\hbar}\int_{\Admun}\psi^{2g}
  \int_{\Admunz}1 \\
= (-)^{g+1}\frac{2}{3}{\binom{2g+3}{2g+2}} P_{3,(3)}(g)L_3(0)\frac{1}{\hbar}.
\end{multline*}

\item[$F_{g_1,g_2}$]
\labellist\small
\pinlabel {$g_1$} [b] at 43 26
\pinlabel {$3$} [b] at 223 13
\pinlabel {$g_2$} [b] at 416 22
\endlabellist
\includegraphics[width=3cm]{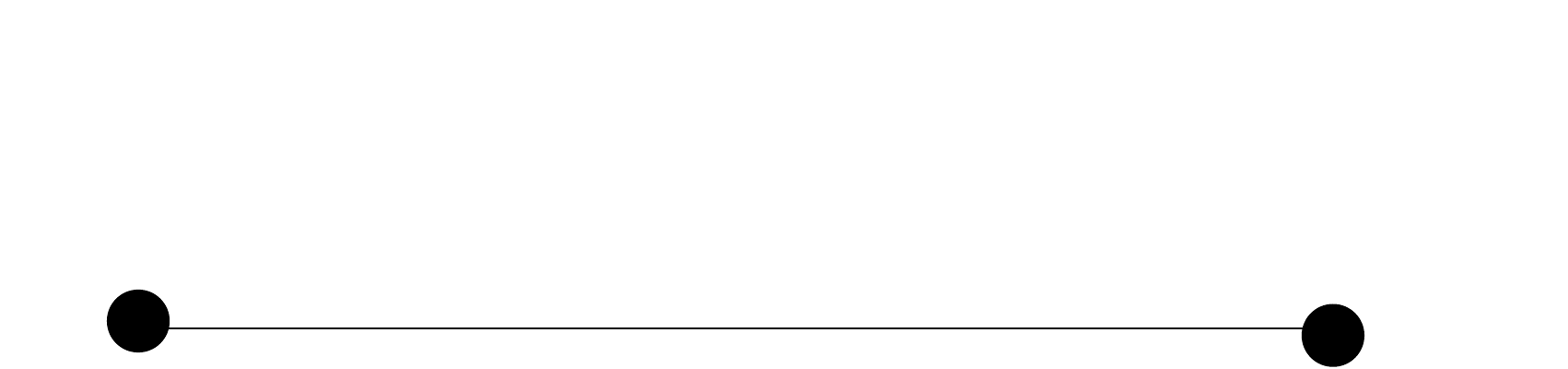}
\begin{align*}
3&{\binom{2g{+}3}{2g_1{+}2}}\int_{\Admungone}
  \frac{(-)^{g_1}\hbar^{2g_1}}{\hbar(\hbar-\psi_3)} \\
&  \hspace{130pt}\int_{\Admungtwo}
  {\frac{\lambda_{g_2}\Lambda_{g_2}(-1)}
  {-\hbar-\psi_3}\left(\frac{2}{9}\hbar^2\right)} \\
&= (-)^{g_1+1}\frac{2}{3}{\binom{2g{+}3}{2g_1{+}2}}
  \frac{1}{\hbar}\int_{\Admungone}\psi^{2g_1} \\
  &\hspace{110pt}\int_{\Admungtwo}\lambda_{g_2}\lambda_{g_2-1}\psi_3
   + \cdots + \lambda_{g_2}\psi_3^{g_2} \\
&\hspace{160pt}= (-)^{g_1+1}\frac{2}{3}{\binom{2g{+}3}{2g_1{+}2}}
  P_{3,(3)}(g_1) L_3(g_2)\frac{1}{\hbar}.
\end{align*}

\item[$F_{0,g}$]
\labellist\small
\pinlabel {$0$} [b] at 43 26
\pinlabel {$3$} [b] at 223 13
\pinlabel {$g$} [b] at 416 22
\endlabellist
\includegraphics[width=3cm]{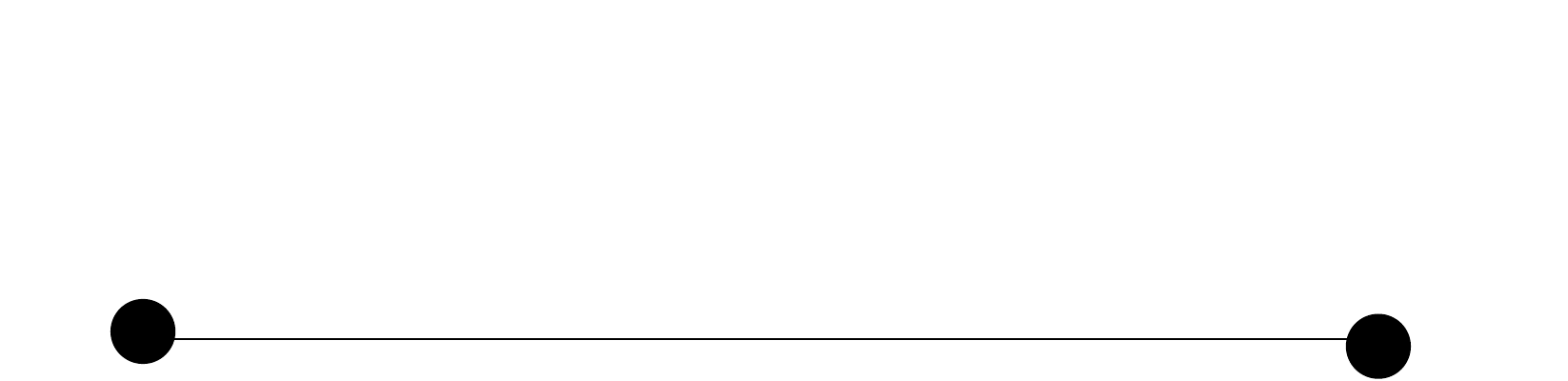}
\begin{align*}
3&{\binom{2g+3}{2}}\int_{\Admunz}
  {\frac{1}{\hbar(\hbar-\psi_3)}} \\
&\hspace{130pt} \int_{\Admun}
  {\frac{\lambda_{g}\Lambda_{g}(-1)}{-\hbar-\psi_3}
  \left(\frac{2}{9}\hbar^2\right)} \\
&= (-)^{g+1}\frac{2}{3}{\binom{2g+3}{2}}
  \frac{1}{\hbar}\int_{\Admunz}1 \\
&\hspace{130pt}\int_{\Admun}\lambda_{g}\lambda_{g-1}\psi_3 + \cdots +
  \lambda_{g}\psi_3^{g} \\
&\hspace{200pt}= - \frac{2}{3}{{2g+3}\choose{2}} P_{3,(3)}(0) L_3(g)\frac{1}{\hbar}.
\end{align*}

\item[$F_{g,\cdot,\cdot}$]
\labellist\small
\pinlabel {$g$} [b] at 51 137
\pinlabel {$2$} [b] at 246 189
\pinlabel {$1$} [b] at 246 65
\endlabellist
\raisebox{-20pt}{\includegraphics[width=3cm]{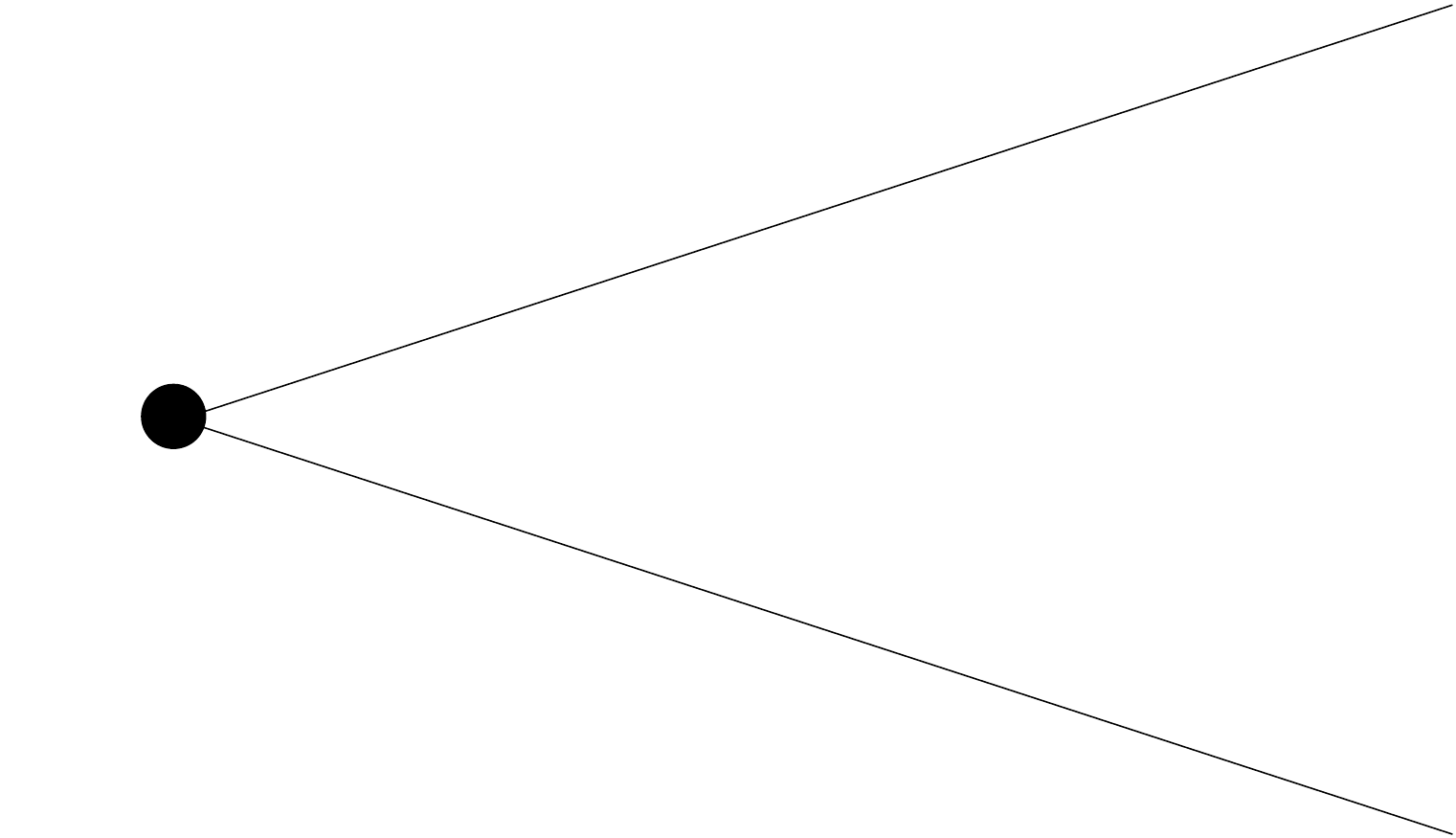}}
\begin{multline*}
{\binom{2g+3}{2g+3}}\int_{\Admunt}
  {\frac{(-)^g \hbar^{2g}}{\hbar(\hbar-\psi_t)}
  \left(\frac{1}{2}\hbar^2\right)} \\
= (-)^{g}\frac{1}{2}\frac{1}{\hbar}\int_{\Admunt}
  \psi_t^{2g+1}
=(-)^g\frac{1}{2}P_{3,(t)}(g)\frac{1}{\hbar}.
\end{multline*}

\item[$F_{g_1,g_2,x}$]
\labellist\small
\pinlabel {$g_1$} [b] at 51 137
\pinlabel {$2$} [b] at 246 189
\pinlabel {$1$} [b] at 246 65
\pinlabel {$g_2$} [b] at 440 262
\endlabellist
\raisebox{-20pt}{\includegraphics[width=3cm]{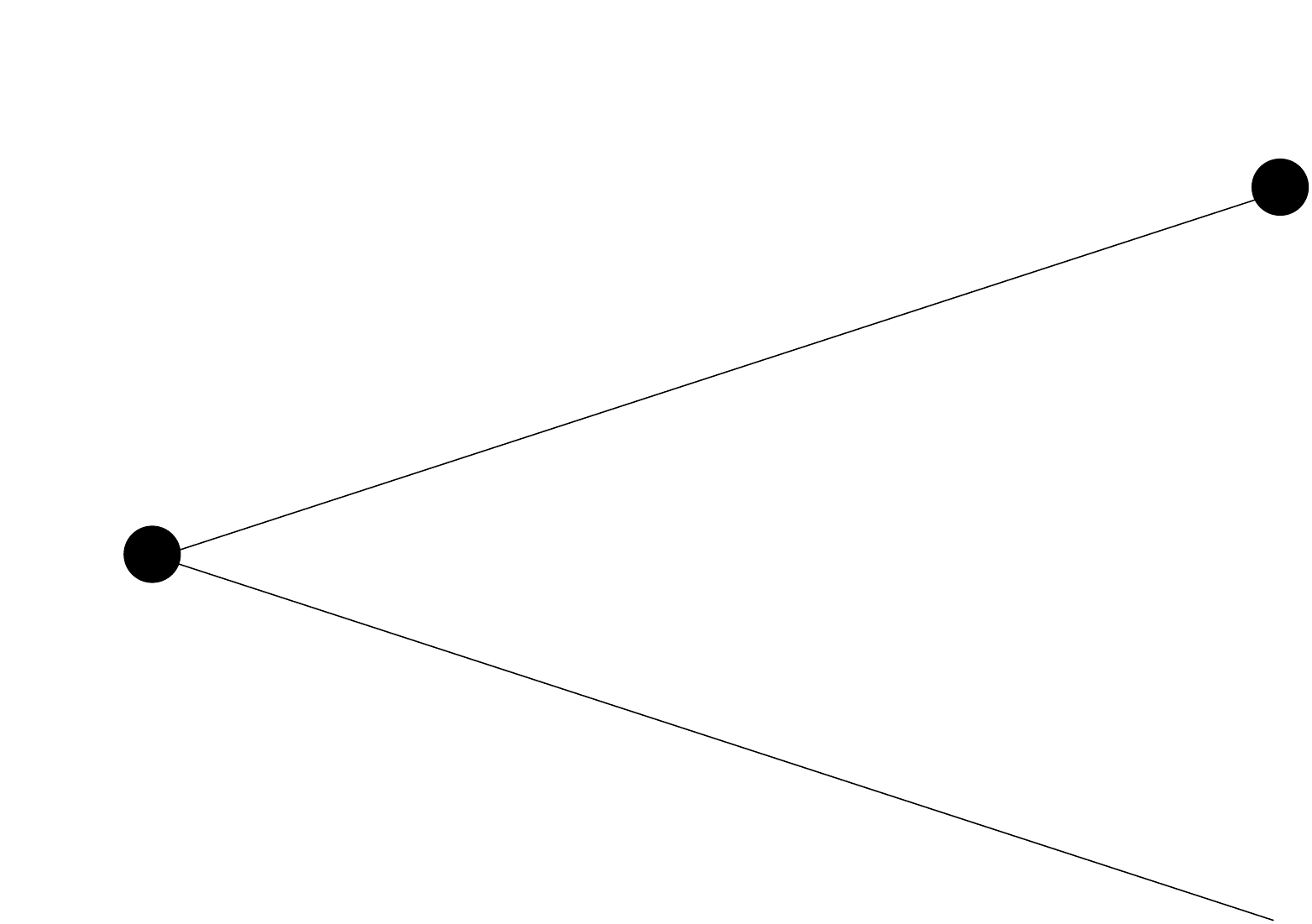}}
\begin{align*}
2&{\binom{2g+3}{2g_1+3}}\int_{\Admuntgone}
  {\frac{(-)^{g_1} \hbar^{2g_1}}{\hbar(\hbar-\psi_t)}} \\
&\hspace{60pt}\int_{\Admuntwogtwo}
  {\frac{\lambda_{g_2}\Lambda_{g_2}(-1)}{-\hbar-\psi_t}
  \left(\frac{\hbar^2}{2}\right)} \\
&\hspace{20pt}= (-)^{g_1}{\binom{2g+3}{2g_1+3}}
  \frac{1}{\hbar}\int_{\Admuntgone}\psi_t^{2g_1+1} \\
&\hspace{130pt}\int_{\Admuntwogtwo}\lambda_{g_2}\lambda_{g_2-1}+\cdots+
  \lambda_{g_2}\psi_t^{g-1} \\
&\hspace{20pt}= (-)^{g_1}{\binom{2g+3}{2g_1+3}}P_{3,(t)}(g_1)L_2(g_2)\frac{1}{\hbar}.
\end{align*}

\item[$F_{0,g,\cdot}$]
\labellist\small
\pinlabel {$0$} [b] at 51 137
\pinlabel {$2$} [b] at 246 189
\pinlabel {$1$} [b] at 246 65
\pinlabel {$g$} [b] at 440 262
\endlabellist
\raisebox{-20pt}{\includegraphics[width=3cm]{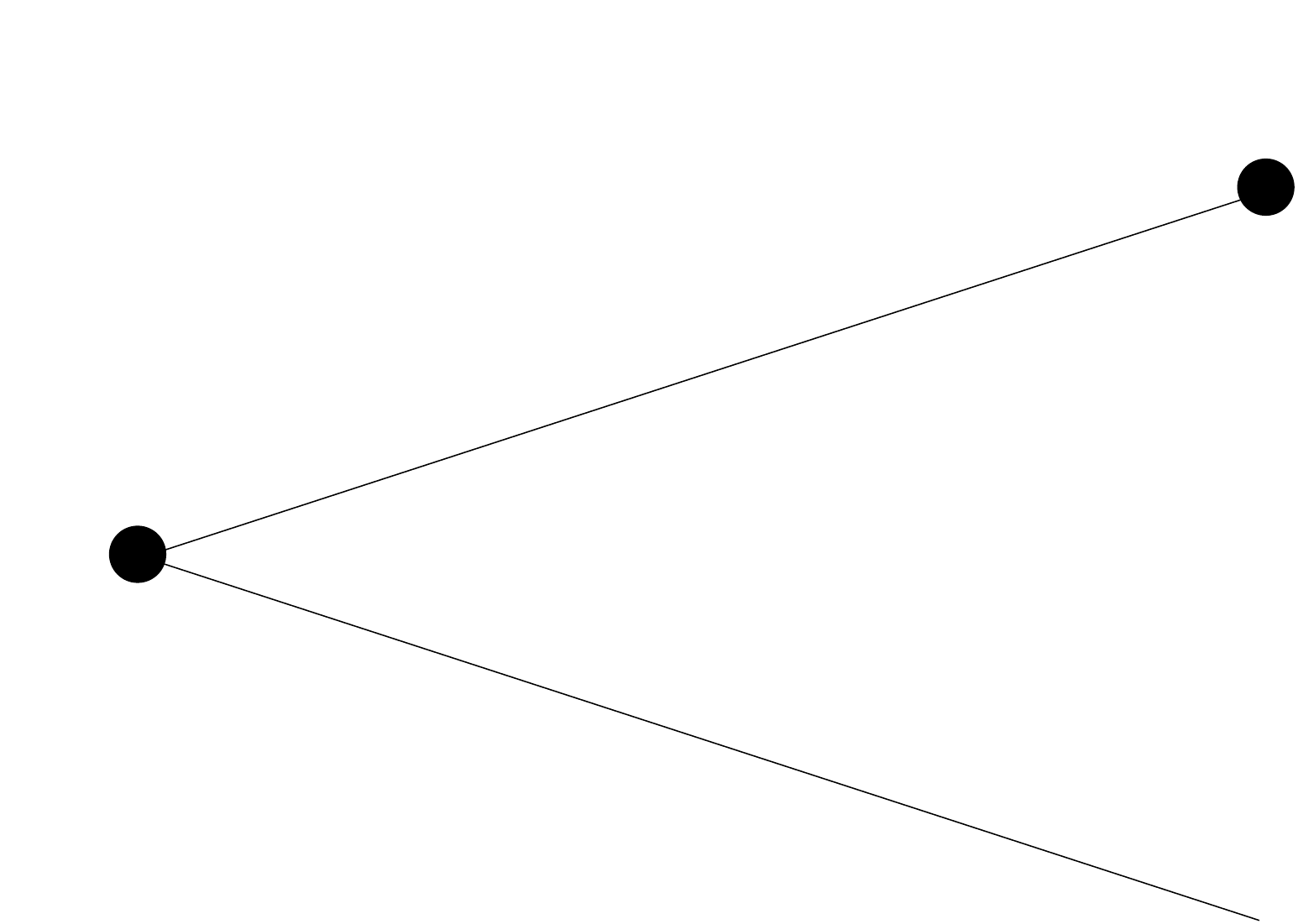}}
\begin{multline*}
2{\binom{2g{+}3}{3}}\int_{\Admuntz}
  {\frac{1}{\hbar(\hbar-\psi_t)}}
  \int_{\Admuntwog}\!\!
  {\frac{\lambda_{g}\Lambda_g(-1)}{-\hbar-\psi_t}
  \left(\frac{1}{2}\hbar^2\right)} \\
={\binom{2g{+}3}{3}}\frac{1}{\hbar}
  \int_{\Admuntz}\psi_t
  \int_{\Admuntwog}\lambda_{g}\lambda_{g-1}+\cdots+
  \lambda_{g}\psi_t^{g-1} \\
= {\binom{2g{+}3}{3}} P_{3,(t)}(0)L_2(g)\frac{1}{\hbar}.
\end{multline*}
\end{description}
Finally, adding everything up, we obtain the following relation:
\begin{multline}
0=\frac{2}{3}\sum_{i=0}^g
  {\binom{2g+3}{2i+1}}(-)^{g-i+1}P_{3,(3)}(g-i)L_3(i)\\[-2ex]
+ \sum_{i=0}^g {\binom{2g+3}{2i}}(-)^{g-i} P_{3,(t)}(g-i)L_2(i).
  \label{lab}
\end{multline}
\subsection{The generating function}
Now for the less deep but more delicate part of our computation: we need
to extract from relation \eqref{lab} a differential equation involving
our desired generating function.

Let's start with a  preliminary  lemma:
\begin{lemma}\label{p}
For all $g\geq0$,
\begin{enumerate} 
\item $ P_{3,(3)}(g)= 3^{2g}.$
\item $ P_{3,(t)}(g)= (3^{2g+2}-1)/2.$
\end{enumerate}
\end{lemma}
\begin{proof}
(1)\qua Consider the map
$$\bfig\morphism<0,-500>[\Admun`\wwbar{M}_{0,2g+3}.;\pi]\efig$$
It's a classical result  that 
$$\int_{\wwbar{M}_{0,2g+3}}\psi_1^{2g}=1.$$ Since our psi class is
just the pull-back of $\psi_1$ on $\wwbar{\mathcal{M}}_{0,2g+3}$, and this space is birational
to its projective coarse moduli space $\wwbar{{M}}_{0,2g+3}$, our lemma is proven if we
show that $\pi$ has degree $3^{2g}$. This is a classic Hurwitz number,
counting the number of degree $3$ covers of the Riemann sphere with a
triple ramification point and simple ramification otherwise.\\ The
problem is purely combinatorial. We are free to choose a three-cycle
in $S_3$ giving the monodromy of the triple point. The triple point
automatically guarantees that our cover is connected. Then we are
free to choose cycles for the first $(2g+1)$ simple ramification
points. The monodromy of the last ramification point is
determined by the fact that the product of all monodromies should be
the identity.  So alltogether we had a choice of $2\cdot3^{2g+1}$
elements of $S_3$. We now need to divide by the conjugation action of
$S_3$ on itself, that geometrically amounts to simply relabelling the sheets of the cover.\\
Finally we obtain the desired $3^{2g}$ non isomorphic covers.

(2)\qua Similarly, we need to count the number of degree $3$ covers of $\proj$ with $2g+4$ simple ramification
points. Paralleling the previous argument, we can choose $(2g+3)$ cycles freely. But we have to beware of 
disconnected covers. These can happen only if we  chose always the same cycle. So in total we have
$3^{2g+3}-3$ choices. Dividing now by $6$ we obtain our claim.
\end{proof}

Let us now  translate  relation \eqref{lab}
in the language of generating functions. Define
\begin{align*}
{\mathcal{L}}_3(x) &:=\sum_{g=0}^\infty\frac{{L}_3(g)}{2g+2!}x^{2g+2}, &
{\mathcal{L}}_2(x) &:=\sum_{g=0}^\infty\frac{{L}_2(g)}{2g+1!}x^{2g+1}, \\
{\mathcal{P}}_{3,(3)}(x)
  &:=\sum_{g=0}^\infty(-)^g\frac{{P}_{3,(3)}(g)}{2g+2!}x^{2g+2}, &
{\mathcal{P}}_{3,(t)}(x)
  &:=\sum_{g=0}^\infty(-)^g\frac{{P}_{3,(t)}(g)}{2g+3!}x^{2g+3}.
\end{align*}
Then our relation \eqref{lab} becomes an ordinary differential equation
on the generating functions:
\begin{equation}
%\begin{imp}
\tfrac{2}{3}{\mathcal{P}}_{3,(3)}{\mathcal{L}}^\prime_3
  -{\mathcal{P}}_{3,(t)}{\mathcal{L}}^\prime_2=0
%\end{imp}
\label{ode3}
\end{equation}
By \fullref{p} we can explicitly describe the Hurwitz numbers'
generating functions
$${\mathcal{P}}_{3,(3)}(x) = \frac{1-\cos(3x)}{9};\qquad
  {\mathcal{P}}_{3,(t)}(x) = \frac{3\sin(x)-\sin(3x)}{6}.$$
Also, we do know the generating function for the degree $2$ theory, hence
$${\mathcal{L}}_2^\prime(x)= \tfrac{d}{dx}\tan{\left(\tfrac{x}{2}\right)}=
  \frac{1}{2\cos^2{\left(\frac{x}{2}\right)}}.$$
Finally, we have reduced our problem to integrating the following:
\begin{equation}
%\begin{imp}
\begin{aligned}
\tilde{\mathcal{L}}^\prime_3(x)&=
\frac{9}{8}\frac{3\sin(x)-\sin(3x)}
  {(1-\cos(3x))\cos^2{\left(\frac{x}{2}\right)}},\\
\tilde{\mathcal{L}}_3(0)&=0.
\end{aligned}
%\end{imp}
\label{explicitode3}
\end{equation}
This ODE integrates to
$${\mathcal{L}}_3(x)=\frac{9}{2}\left(\frac{1}{4\cos^2{\left(\frac{x}{2}\right)}-1}-\frac{1}{3}\right).$$
Now let us remember that the generating function $\mathcal{I}_3(x)$ is smply  $(2/9)\mathcal{L}_3(x)$. After just a little bit of trigonometry clean-up we obtain:
\begin{equation}
%\begin{vimp}
\mathcal{I}_3(x)=
  \frac{4}{3}\frac{\sin^3{\left(\frac{x}{2}\right)}}
  {\sin{\left(\frac{3x}{2}\right)}}
\label{deg3}
%\end{vimp}
\end{equation}
\subsection{A corollary}
In a completely similar fashion to degree $2$, it is possible to obtain from the previous computation the generating function for the integrals:
$$J_3(g)=\int_{\Admtinf} c_{2g+4}(\F).$$
The answer is
\begin{equation}
%\begin{vimp}
\mathcal{J}_3(x)= \sum_{g=0}^\infty
  \left(\frac{J_3(g)}{2g+4!}\right)x^{2g+4}= 3\mathcal{I}(x)^2=
  \frac{16}{3}\frac{\sin^6{\bigl(\frac{x}{2}\bigr)}}
  {\sin^2{\bigl(\frac{3x}{2}\bigr)}}.
%\end{vimp}
\end{equation}

\bibliographystyle{gtart}
\bibliography{link}

\end{document}